\documentclass[nonblindrev]{informs3} % current default for manuscript submission

\OneAndAHalfSpacedXI
\renewcommand{\arraystretch}{1.25}
\usepackage{setspace}
\usepackage{natbib}
\bibpunct[, ]{(}{)}{,}{a}{}{,}%
\def\bibfont{\small}%

{}
{}
{}
\usepackage{lmodern}
\usepackage{booktabs}
\TheoremsNumberedThrough     % Preferred (Theorem 1, Lemma 1, Theorem 2)
\ECRepeatTheorems

\EquationsNumberedThrough    % Default: (1), (2), ...
\usepackage[colorlinks=true, allcolors=black]{hyperref}
\usepackage{makecell}
\usepackage{bm}

\defcitealias{Bagnell2001SolvingUM}{Bagnell et al. (2001)}
\defcitealias{wiesemann2013robust}{Wiesemann et al. (2013)}
\defcitealias{mannor2016robust}{Mannor et al. (2016)}
\defcitealias{gohDataUncertaintyMarkov2018}{Goh et al. (2018)}
\defcitealias{delimpaltadakis2023interval}{Delimpaltadakis et al. (2023)}
\defcitealias{black2022robust}{Black et al. (2023)}
\defcitealias{song2023decision}{Song et al. (2024)}
\defcitealias{chen2019distributionally}{Chen et al. (2019)}
\defcitealias{chen2020robust}{Chen et al., 2020}
\defcitealias{weissman2003inequalities}{Wiesemann et al., 2003}
\defcitealias{ho2018fast}{Ho et al. (2018)}
\defcitealias{ho2021partial}{Ho et al. (2021)}
\defcitealias{behzadian2021fast}{Behzadian et al. (2021)}
\defcitealias{russel2019optimizing}{Russel et al. (2019)}
\defcitealias{klabjan2013robust}{Klabjan et al. (2013)}
\defcitealias{wang2023finite}{Wang et al. (2023)}
\defcitealias{ho2022robust}{Ho et al. (2022)}
\defcitealias{mannor2012lightning}{Mannor et al. (2012)}
\defcitealias{ben2009robust}{Ben et al., 2009}
\defcitealias{tirinzoni2018policy}{Tirinzoni et al. (2018)}
\defcitealias{wang2023policy}{Wang et al. (2018)}
\defcitealias{long2023robust}{Long et al., 2023}
\defcitealias{ghavamzadeh2016safe}{Ghavamzadeh et al. (2016)}
\defcitealias{laroche2019safe}{Laroche et al. (2019)}
\defcitealias{simao2020safe}{Simao et al. (2020)}
\defcitealias{lim2013reinforcement}{Lim et al. (2013)}
\defcitealias{ruan2023robust}{Ruan et al. (2013)}
\defcitealias{auer2009near}{Auer et al., 2009}
\defcitealias{derman2020bayesian}{Derman et al. (2020)}
\defcitealias{o2018uncertainty}{O'Donoghue et al. (2018)}
\defcitealias{wang2023robust}{Wang et al. (2023)}
\defcitealias{kumar2024policy}{Kumar et al. (2023)}
\defcitealias{huang2021convergence}{Huang et al., 2021}
\defcitealias{zhang2023regularized}{Zhang et al. (2023)}
\defcitealias{derman2021twice}{Derman et al. (2021)}
\defcitealias{li2022quantile}{Li et al. (2022)}
\defcitealias{mckenna2020approximate}{McKenna et al., 2020}
\defcitealias{song1993inventory}{Song \& Zipkin, 1993}
\defcitealias{che2023deep}{Che et al., 2023}
\defcitealias{fan2022optimal}{Fan et al., 2022}
\defcitealias{bushaj2023simulation}{Bushaj et al., 2023}
\defcitealias{kosanoglu2022deep}{Kosanoglu et al., 2022}
\defcitealias{matsuo2022deep}{Matsuo et al., 2022}
\defcitealias{rolf2022review}{Rolf et al., 2023}
\defcitealias{zuccotto2024reinforcement}{Zuccotto et al., 2024}
\defcitealias{ben2002robust}{Ben-Tal \& Nemirovski, 2002}
\defcitealias{bertsimas2004price}{Bertsimas \& Sim, 2004}
\defcitealias{goerigk2016algorithm}{Goerigk \& Schobel, 2016}
\defcitealias{scarf1957min}{Scarf et al., 1957}
\defcitealias{givan2000bounded}{Givan et al., 2000}
\defcitealias{white1986parameter}{White III \& ElDeib, 1986}
\defcitealias{white1994markov}{White III \& ElDeib, 1994}
\defcitealias{mannor2007bias}{Mannor et al. (2007)}
\defcitealias{agussurja2019state}{Agussurja et al., 2019}
\defcitealias{shi2021timing}{Shi et al., 2021}
\defcitealias{drent2024condition}{Drent et al., 2021}
\defcitealias{cowan2018reinforcement}{Cowan et al., 2018}
\defcitealias{croonenborghs2007online}{Croonenborghs et al., 2007}
\defcitealias{badrinath2021robust}{Badrinath \& Kalathil, 2021}
\defcitealias{hans2008safe}{Hans et al, 2008}
\defcitealias{polo2011safe}{Polo \& Rebollo, 2011}
\defcitealias{wachi2020safe}{Wachi \& Sui, 2020}
\defcitealias{farahmand2011regularization}{Farahmand, 2011}
\defcitealias{zhang2018study}{Zhang et al., 2018}
\defcitealias{badings2023decision}{Badings et al., 2023}
\defcitealias{keith2021survey}{Keith \& Ahner, 2021}
\defcitealias{mannor2019data}{Mannor \& Xu, 2019}
\defcitealias{philpott2008convergence}{Philpott \& Guan, 2008}
\defcitealias{zou2019stochastic}{Zou et al., 2019}
\defcitealias{georghiou2019robust}{Georghiou et al., 2019}
\defcitealias{daryalal2023two}{Daryalal et al., 2023}
\defcitealias{daryalal2024lagrangian}{Daryalal et al., 2024}
\defcitealias{bertsimas2019adaptive}{Bertsimas et al., 2019}
\defcitealias{delage2010distributionally}{Delage \& Ye, 2010}
\defcitealias{ninh2021robust}{Ninh, 2021}
\defcitealias{wiesemann2014distributionally}{Wiesemann et al., 2014}
\defcitealias{yu2015distributionally}{Yu and Xu (2015)}
\defcitealias{pmlr-v130-behzadian21a}{Behzadian et al., 2021}
\defcitealias{blanchet2015unbiased}{Blanchet \& Glynn, 2015}
\defcitealias{blanchet2019unbiased}{Blanchet et al., 2019}
\defcitealias{shi2024curious}{Shi et al. (2023)}
\defcitealias{ben2004adjustable}{Ben-Tal et al., 2004}
\defcitealias{yanikouglu2019survey}{Yan{\i}ko{\u{g}}lu et al., 2019}
\begin{document}
%%%%%%%%%%%%%%%%

% Outcomment only when entries are known. Otherwise, leave it as is and
%\setcounter{page}{1}
%\DOI{10.1287/xxxx.0000.0000}%


\RUNTITLE{Wenfan Ou and Sheng Bi}

% \TITLE{Bundling Information Goods of Decreasing Value}
\TITLE{Sequential Decision-Making under Uncertainty: A Robust MDPs review}

% NOTE: Authors with the same affiliation, if the order of authors allows,
\ARTICLEAUTHORS{%
\AUTHOR{Wenfan Ou}
\AFF{School of Information Management and Engineering, Shanghai University of Finance and Economics, Shanghai 200437, China,  \EMAIL{ouwenfan@stu.sufe.edu.cn}} %, \URL{}}
\AUTHOR{Sheng Bi}
\AFF{School of Information Management and Engineering, Shanghai University of Finance and Economics, Shanghai 200437, China, \EMAIL{bisheng@sufe.edu.cn}}
} % end of the block

\ABSTRACT{%
Fueled by advances in both robust optimization theory and reinforcement learning (RL), robust Markov Decision Processes (RMDPs) have garnered increasing attention due to their powerful capability for sequential decision-making under uncertainty. In this paper, we provide a comprehensive overview of the theoretical foundations and recent developments in RMDPs, with a particular emphasis on ambiguity modeling. We examine the ``rectangular assumption", a key condition ensuring computational tractability in RMDPs but often resulting in overly conservative policies. Three widely used rectangular forms are summarized, and a novel proof is provided for the NP-hardness of non-rectangular RMDPs. We categorize RMDP formulation approaches into parametric, moment-based, and discrepancy-based models, analyzing the trade-offs associated with each representation. Beyond the traditional scope of RMDPs, we also explore recent efforts to relax rectangular assumptions and highlight emerging trends within the RMDP research community. These developments contribute to more practical and flexible modeling frameworks, complementing the classical RMDP results. Relaxing rectangular assumptions tailored to operations management is a promising area for future research, and there are also opportunities for further advances in developing fast algorithms and provably robust RL algorithms.
% Enter your abstract
}%

%\KEYWORDS{deterministic inventory theory; infinite linear programming duality;
%  existence of optimal policies; semi-Markov decision process; cyclic schedule}


% \HISTORY{This paper was
% first submitted on April 12, 1922 and has been with the authors for
% 83 years for 65 revisions.}

\maketitle
%%%%%%%%%%%%%%%%%%%%%%%%%%%%%%%%%%%%%%%%%%%%%%%%%%%%%%%%%%%%%%%%%%%%%%

% Samples of sectioning (and labeling) in MNSC
% NOTE: (1) \section and \subsection do NOT end with a period
%       (2) \subsubsection and lower need end punctuation
%       (3) capitalization is as shown (title style).
%
%\section{Introduction.}\label{intro} %%1.
%\subsection{Duality and the Classical EOQ Problem.}\label{class-EOQ} %% 1.1.
%\subsection{Outline.}\label{outline1} %% 1.2.
%\subsubsection{Cyclic Schedules for the General Deterministic SMDP.}
%  \label{cyclic-schedules} %% 1.2.1
%\section{Problem Description.}\label{problemdescription} %% 2.

% Text of your paper here

\section{Introduction}\label{sec1}

Sequential decision-making under uncertainty is common in almost every scientific domain and practical activity. Typically, a decision-maker selects an \textit{action} based on the available information, referred to as a \textit{state}, at each decision epoch, and then obtains uncertain outcomes and a new state for the next epoch. This feedback loop allows multistage problems to be modeled as a Markov Decision Process (MDP), provided the Markov property holds, where the next state and reward depend only on the current state and action. Different from the methods that aim to capture exploration and exploitation trade-offs (e.g., multi-armed bandits), MDPs focus on \textit{planning}, often assuming the reward function and transition probabilities are known \citep{puterman2014markov}. Given probabilistic knowledge of the uncertainty, the decision-maker can employ dynamic programming (DP) or linear programming (LP) techniques to determine the optimal policy of the multistage problem. These methods are typically more computationally efficient than brute-force enumeration, such as decision trees. Due to both the broad applicability of the \textit{state} concept and the inherent recursive structure of the Bellman equation (as detailed later), numerous problems across various fields can be formulated as MDPs, such as economics \citep{parkes2003mdp, katehakis2012optimal}, transportation \citepalias{agussurja2019state, mckenna2020approximate}, inventory control \citep{song1993inventory, feinberg2022optimality}, healthcare \citepalias{shi2021timing, fan2022optimal}, and manufacturing \citepalias{drent2024condition}. These advantages enable a concise representation of system dynamics (through state transitions) while also facilitating problem decomposition and efficient computation.

In this survey, we focus on the problems where the transition probabilities are unknown. We assume rewards either entirely depend on the current state and action or are governed by the transition probabilities. Full knowledge of the transition probabilities is usually unavailable in the real world. Under a data-driven paradigm, the decision-maker can only estimate model parameters or distributions based on historical data, which introduces additional \textit{external uncertainty} (distinguished from the \textit{internal uncertainty} due to the stochastic nature of MDPs). With the deviated estimation or belief affected by external uncertainty, the decision-maker should anticipate disappointment on average, owing to the optimization-based selection process. This is the so-called the \textit{optimizer's curse} \citep{smith2006optimizer}. Even worse, since sequential decision-making can amplify decision errors across stages, inaccurate estimations often result in significantly worse outcomes than in single-stage problems that involve only one decision point. We refer to the case study in \citetalias{mannor2007bias} to demonstrate the influence of external uncertainty under the optimization process.

One effective approach to mitigate the impact of external uncertainty is robust optimization (RO) \citepalias{ben2002robust, goerigk2016algorithm}. Traditional RO typically involves single-period problems by seeking a solution that is optimal under the worst case (constrained by certain conditions), ensuring performance despite incomplete probabilistic information \citep{scarf1957min, bertsimas2004price}. An important extension to multi-stage settings is adjustable robust optimization (ARO), which reduces conservatism by allowing decision-makers to adapt to new information \citepalias{yanikouglu2019survey}. However, its general formulation is proven computationally intractable \citepalias{ben2004adjustable}, and solving ARO problems necessitates exploiting problem-specific structures or employing approximation methods.

Regarding the computationally tractable recursive structure of MDP, another natural extension to multi-stage settings is \textit{robust MDP} (RMDP). In recent years, driven by breakthroughs in theory and practical needs of data-driven paradigms, RO has seen substantial progress, and RMDP has also been emerging within operations research and computer science communities. Unlike conventional MDPs, RMDPs are tailored for circumstances where only partial information about the true underlying parameters or probability distributions is available, rather than complete knowledge. The decision-maker can construct a collection of possible parameter values (resp. distributions), called an \textit{uncertainty set} (resp. \textit{ambiguity set}), based on the belief that the partial information is credible. This allows the decision-maker to determine the optimal solution under the worst-case scenario over the uncertainty or ambiguity set.

The investigation of RMDPs can be traced back to the 1970s \citep{satia1973markovian} when they were known as MDPs with imprecisely known parameters (MDPIPs) \citepalias{white1986parameter, white1994markov, givan2000bounded}. MDPIPs assume the transition probabilities are constrained by a finite number of linear inequalities that form a linear program (e.g., each transition probability is bounded within a reasonable constant range). This is a simple and natural way to characterize uncertainty, and the polytopic models are computationally tractable and can be solved by LP techniques. However, this approach also brings some drawbacks. Polytopic models not only present computational challenges due to the equality constraint imposed on (conditional) transition probabilities but also result in overly conservative policies stemming from the statistically poor representations of uncertainty. For instance, a simplistic uncertainty set for an MDPIP might constrain each (conditional) transition probability to be limited in range $[0.1,0.9]$, with their summation equaling $1$. In this uncertainty set, each (conditional) transition probability is considered individually, except for the summation constraint. However, as previously mentioned, additional information about the overall transition probabilities can be derived from historical data or empirical probabilities beyond these ranges, such as their shape or symmetry.

To address this issue, seminal papers by \cite{iyengarRobustDynamicProgramming2005} and \cite{nilimRobustControlMarkov2005} independently show that a robust formulation with the $(s, a)$-rectangular assumption can be solved recursively via the robust Bellman equation (details provided later), thereby extending the results from non-robust DP theory. The rectangularity assumption allows the problem to be decomposed into independent subproblems. Therefore, they can adapt the value iteration and policy iteration for MDPs into \textit{robust} counterparts, enabling efficient computation. This formulation is also referred to as \textit{robust} DP and serves as a canonical form of RMDPs thereafter. Both papers introduce two imperative concepts. Firstly, they introduce uncertainty sets with statistical metrics, especially $\phi$-divergences, norms, and likelihood region ambiguity sets, which outperform polytopic models. Secondly, these works coincidentally emphasize that the rectangularity assumption is essential for the tractability of RMDPs, albeit without providing concrete proof. Bridging this theoretical gap, \citetalias{Bagnell2001SolvingUM} and \citetalias{wiesemann2013robust} independently establish the NP-hardness of general (non-rectangular) RMDPs through distinct approaches. 

While RMDPs have utilized statistical metrics to construct uncertainty sets and enhance representations of uncertainty, these methods do not incorporate prior distribution information about the uncertainty. Specifically, whether using an entropy metric or a likelihood region, the obtained uncertainty set only contains a collection of possible parameter realizations (i.e., each conditional transition probability) such that only their support is known. This is why such a collection is termed an \textit{uncertainty set}. However, prior information is often available in practical applications, such as domain knowledge. Considering the potential conservatism of traditional robust approach \citep{thiele2010note, delage2010percentile} and advances in distributionally robust optimization (DRO) \citep{popescu2007robust, delage2010distributionally, goh2010distributionally}, \cite{xu2012distributionally} first propose the framework of distributionally robust MDPs (DRMDPs), where the optimal decision is against the worst distribution among a collection termed as an \textit{ambiguity set}.  Compared to traditional RMDPs, the new framework more conveniently integrates statistical information to mitigate conservatism. Both \cite{xu2012distributionally} and \citetalias{yu2015distributionally} demonstrate that DRMDPs can be reduced to expected standard RMDPs under certain assumptions, thus retaining tractability. These research findings have laid a solid theoretical foundation for RMDPs and enable the thriving development of the RMDP community.

Besides the exploration in the stochastic optimization area, the literature on RMDPs has been significantly enriched by diverse fields recently. Researchers across different fields combine RMDPs with other techniques from various perspectives, including online learning \citepalias{croonenborghs2007online, badrinath2021robust,cowan2018reinforcement}, safe learning \citepalias{hans2008safe, polo2011safe, wachi2020safe}, dynamic risk measure \citep{bauerle2014more, yu2022risk, rockafellar2000optimization}, regularization \citepalias{farahmand2011regularization, zhang2018study}, etc. Despite the surge in the literature on RMDPs in recent years, a systematic review of the formulation of RMDPs remains scarce, with only a handful of tutorials or general sequential decision-making under uncertainty reviews available \citepalias{mannor2019data, keith2021survey, badings2023decision}. We believe that the absence of a specialized review hinders the development of RMDP theory and applications. To bridge this gap, we provide an up-to-date review of the RMDP/DRMDP formulations with uncertain transitions. This is a common class of RMDP problems that have been widely considered in research and applications within the management science and operations research community. By surveying significant results and studies in this domain, we aim to establish a solid foundation for future research and practical implementations in RMDPs.

The remainder of the paper is organized as follows. In Section \ref{sec2}, we briefly introduce the preliminary of RMDPs. In particular, we summarize the definitions of rectangularity and provide a new proof for the NP-hardness of non-rectangular RMDPs. In Sections \ref{sec3}, \ref{sec4}, and \ref{sec5}, we introduce, respectively, three popular types of RMDPs: parametric RMDPs, moment-based RMDPs, and discrepancy-based RMDPs. In these sections, we focus on generic formulations and key theoretical advances. In Section \ref{sec6}, we investigate methodologies for modeling coupled uncertainty that violates rectangularity assumptions.  Finally, we review some novel works in Section \ref{sec7} that have sparked our interest beyond the traditional mini-max framework and conclude in Section 8.

\section{Preliminaries}\label{sec2}

We denote by $[T]=\{1,2,..., T\}$ the set of positive running indices up to $T$, and use $\bigotimes$ to represent the Cartesian product. Let $\Xi:=\bigotimes_{t\in[T]}\Xi_t$ be the entire sample space, where $\Xi_t$ is the sample space at stage $t$. Let $\mathcal{F}$ be the $\sigma$-algebra of $\Xi$, and $\mathcal{F}$ consists of all subsets of $\Xi$. Exactly,  the $\sigma$-algebra $\mathcal{F}$ is a set of events, where an event can be viewed as a scenario trajectory in the context of this paper. For the events that can be observed up to $t$, we denote the set of them by the filtration $\mathcal{F}_t$. Thus, we have $\left\{ \emptyset,\Xi\right\} =\mathcal{F}_{0}\subseteq\mathcal{F}_{1}\subseteq\cdots\subseteq\mathcal{F}_{T}=\mathcal{F}$.  One can think of this information structure also in terms of a scenario tree, $\mathcal{F}_t$ is generated by the partitions corresponding to nodes at stage $t$. Let $\mathfrak{M}(\Xi,\mathcal{F})$ be a set of probability measures on measurable space $(\Xi,\mathcal{F})$. We denote  $\mathcal{P}\subseteq\mathfrak{M}(\Xi,\mathcal{F})$ by the set of probability measures we are interested in, and $\mathbb{P}\in \mathcal{P}$ is a probability measure on $(\Xi,\mathcal{F})$. Similarly, we specify $\mathcal{P}_{t}=\left\{ \mathbb{P}(\cdot|\mathcal{F}_{t})\ |\ \mathbb{P}\in\mathcal{P}\right\}$ where $\mathbb{P}(\cdot|\mathcal{F}_t)$ is the probability measure $\mathbb{P}$ conditioned on the filtration $\mathcal{F}_t$. In this survey, we particularly focus on the most common case that the (conditional) transition probabilities are time-invariant. This assumption is common in the literature, such as inventory and queueing systems, implying that system dynamics do not change over time. Except for otherwise noted, we assume stationary environments thereafter. As  $\Xi$ is finite, without loss of generality, we can appropriately derive the corresponding probabilities $\mathbf{p}$, a (discrete) distribution $\mu$ or a probability measure $P$ on a single stage measurable space $(\Xi_t,\mathcal{F}^t)$ for all $t\in[T]$ with respect to $\mathbb{P}$, where $\mathcal{F}^t$ is the corresponding $\sigma$-algebra with respect to $\Xi_t$. Let $ \Delta^{d}=\left\{ \mathbf{p}\ |\ p_{1}+\cdots+p_{d}=1,\ p_{i}\ge0\ \forall i\in[d]\right\} \in\mathbb{R}^{d}$ be the probability simplex, which represents the set of all valid probability distributions over $d$ outcomes, and $\boldsymbol{1}(\cdot)$ be an indicator function. We use superscript $0$ to denote true underlying true value, e.g., $\mathbf{p}^0$, and $\hat{\cdot}$ denotes the estimated or nominal one, e.g., $\hat{\mathbf{p}}$.

As discussed in the introduction, we focus on sequential decision-making problems where the decision-maker observes the current state of the system, implements a feasible action, and stochastically transits to a new state with associated rewards. We consider a discounted case where the planning horizon $T$ can be infinite. Particularly, the formulations with $T=\infty$ often lead to the existence of stationary policies, which simplifies analysis and contributes to comprehending long-run behavior. When $T=\infty$, we usually require a discount factor $\gamma<1$ to ensure that the interested problems are well-defined. Let $\mathcal{S}$ denote the state space, and $\mathcal{A}$ denote the action space. We use a subscript $s$ to denote the value associated with the state $s$, e.g., $\mathcal{A}_s$. Here, we assume the state space and action space at each stage are time-invariant unless otherwise specified, i.e., $\Xi_t=\mathcal{S}\times\mathcal{A}$ for all $t\in[T]$. For $s\in \mathcal{S}$, $a\in\mathcal{A}$, let  $p_{sas^{\prime}}:=P(s^{\prime}|s,a)$ denote the probability of transitioning to the next state $s^{\prime}\in \mathcal{S}$. These probabilities form the transition probability matrix $\mathbf{p}\in(\Delta^{|\mathcal{S}|})^{|\mathcal{S}|\times|\mathcal{A}|}$. When the transition $(s, a,s')$ occurs, the decision-maker receives an immediate reward $R(s, a, s')$. For simplicity, we assume $R:\mathcal{S}\times\mathcal{A}\times\mathcal{S}\rightarrow \mathbb{R}$ is a known and deterministic mapping except for additional notes in the context. Given the discount factor $0<\gamma\le1$, such a multistage problem can be modeled as a MDP, which is represented as a six-element tuple $<T,\mathcal{S},\mathcal{A}, R,\mathbf{p},\gamma>$. The goal of decision-makers is to maximize the cumulative discounted rewards by optimizing a policy $\pi$:
\begin{equation}
    \label{multistageobj}
    \max_{\pi}\mathbb{E}_{\mathbf{p}}^{\pi}\left[{\textstyle \sum\limits _{t=1}^{T}\gamma^{t-1}R(s_{t},a_{t},s_{t+1})+\gamma^{T}R(s_{T+1})}\right]
\end{equation}
where $R(s_{T+1})$ is the final reward function. A policy is a mapping from states to actions (deterministic) or a probability distribution over actions (stochastic). For example, given the state $s\in \mathcal{S}$, a deterministic policy $\pi(s)=a$ with probability 1, whereas a stochastic policy takes action $a$ with probability $\pi(a\mid s)$.

Once at a given state, we denote the value-to-go function $V_t(\cdot):\mathcal{S}\rightarrow \mathbb{R}$ as the mapping that quantifies the expected future reward generated by policy $\pi$ at stage $t$: 
\begin{equation}
	\label{valuefunction}
	V_t(s)=\max_{\pi}\mathbb{E}_{\mathbf{p}}^{\pi}\left[\sum_{\tau=t}^T \gamma^{\tau-1}R(s_{\tau}, a_{\tau}, s_{\tau+1})+\gamma^TR(s_{T+1})\,\middle|\, s_t=s\right].
\end{equation} Fortunately, by leveraging the value-to-go function as given in (\ref{valuefunction}), solving (\ref{multistageobj}) is equivalent to recursively solving the following equations, known as the Bellman equations:
\begin{equation}
    \label{bellmaneuqation}
    V_t(s)=\max_{\pi} \mathbb{E}_{\mathbf{p}}^{\pi}\left[R(s_t, a_t, s_{t+1})+\gamma V_{t+1}(s_{t+1})\mid s_t=s\right], \ \forall\ t\in[T]
\end{equation}
where $V_{T+1}$ is the terminal value which takes the value $-R(s_{T+1})$ conventionally. 

The Bellman equations decompose the multistage optimization problem into smaller subproblems by recursively relating the value of a state to the values of successor states. The advantages of solving the Bellman equation as in (\ref{bellmaneuqation}) are two-fold. First, the optimal equation retains the same form as the (optimal) Bellman equation, is that
\begin{equation}
	V^{*}(s)=\max_{\pi} \mathbb{E}_{\mathbf{p}}^{\pi}\left[R(s, \pi(s), s')+\gamma V^*(s')\right]
\end{equation} Second, its recursive formulation brings significant computational benefits. For the finite horizon case, straightforwardly applying DP techniques, we can compute the optimal value via backward induction. For instance, in inventory management, backward induction can optimize ordering policies by considering expected immediate newsvendor costs and the value-to-go function at the next stage. For the infinite horizon case, we can compute the optimal value and policy via iterative methods, e.g., value iteration (VI) or policy iteration (PI) \citep{howard1960dynamic}. For example, we operate the maximization as on the right-hand side of (\ref{bellmaneuqation}) in VI and then update the value for each state. This process is repeated until the maximum deviation between two consecutive iterations is smaller than a predefined tolerance level $\epsilon$. Although the ``curse of dimensionality" \citep{bellman1966dynamic} is a known challenge, small to medium-scale MDPs can still be solved efficiently. For more detailed fundamentals, we refer to \cite{puterman2014markov}.

Note that there are various representations of MDPs in different contexts, such as assuming known initial state distributions or continuous state spaces. The above definition of MDP captures the core components and remains widely applicable across this survey. Additional extensions will build on this formulation if necessary.

\begin{remark}
\label{probability_remark}
    The transition probability matrix $\mathbf{p}$ defined above is associated with the transition tuple $(s, a, s')$, but plenty of literature assumes exogenous randomness. For example, in inventory management problems, the demand is often not affected by inventory position and order quantity. In this case, $\mathbf{p}$ may be independent of states and actions and it can be represented by a distribution directly.
\end{remark}

\subsection{(Distributionally) Robust Markov Decisions Processes}

In this survey, we restrict our attention to MDPs with uncertain transition probabilities. To ensure generality, we consider the probability space $(\Xi,\mathcal{F})$, where each $\Xi_t=\mathcal{S}\times\mathcal{A}$, and the set $\mathcal{P}$ consists of probability measures of interest. However, to align with the literature and for simplicity, we also slightly abuse the notation $\mathcal{P}$ to represent the set of corresponding transition probabilities $\mathbf{p}$ or distribution $\mu$, which means $\mathcal{P}$ denotes the collection of uncertainty in this survey, including uncertainty set and ambiguity set. Formally, we consider RMDPs and DRMDPs, defined as follows:

\begin{definition}[Robust MDPs]
    An RMDP is defined as a tuple $<T,\mathcal{S},\mathcal{A},\mathcal{P},R,\gamma>$, where the transition probability matrix $\mathbf{p}\in (\Delta^{|\mathcal{S}|})^{|\mathcal{S}|\times|\mathcal{A}|}$ is unknown but belongs to a known collection $\mathcal{P}$ of parameters, called uncertainty set. More precisely, RMDPs can be expressed as 
    \begin{equation}
    \label{rmdp}
    \max_{\pi}\min_{\mathbf{p}\in\mathcal{P}}\mathbb{E}^{\pi}_{\mathbf{p}}\left[{\textstyle \sum\limits _{t=1}^{T}\gamma^{t-1}R(s_{t},a_{t},s_{t+1})+\gamma^{T}R(s_{T+1})}\right].
    \end{equation}
\end{definition}

\begin{definition}[Distributionally robust MDPs]
    An DRMDP can be represented by a tuple $<T,\mathcal{S},\mathcal{A},\mathcal{P},R,\gamma>$, where the transition probability matrix $\mathbf{p}$ follows a distribution $\mu$, denoted by $\mathbf{p}\sim \mu$. The distribution $\mu$ is unknown but belongs to a known collection $\mathcal{P}$ of probability distributions, called the ambiguity set. More specifically, DRMDPs can be expressed as
    \begin{equation}
    \max_{\pi}\min_{\mu\in\mathcal{P}}\mathbb{E}^{\pi}_{\mathbf{p}\sim \mu}\left[{\textstyle \sum\limits _{t=1}^{T}\gamma^{t-1}R(s_{t},a_{t},s_{t+1})+\gamma^{T}R(s_{T+1})}\right].
    \end{equation}
\end{definition}

While both RMDPs and DRMDPs consider unknown transition probabilities, they handle this uncertainty in different ways. In RMDPs, the transition probability matrix $\mathbf{p}$ is regarded as an uncertain ``parameter", analogous to the uncertain outcomes in RO. The corresponding uncertainty set $\mathcal{P}$ contains all possible realizations of $\mathbf{p}$. The goal of (\ref{rmdp}) is to derive the optimal policy $\pi^*$ against the worst realization. Essentially, the only information available about the uncertainty in RMDPs is the support of $\mathbf{p}$, although this support can be refined using various methods, e.g., statistical metrics. In contrast, DRMDPs assume that additional prior distribution information about $\mathbf{p}$ is available. For instance, one might have a belief that an arbitrary transition probability $p_{sas'}$ does not exceed 0.8 with at least 90\% probability. Such belief can be readily incorporated into the ambiguity sets of DRMDPs, but not into the uncertainty sets of RMDPs. This distinction about how to treat transition probabilities is noteworthy, bringing various perspectives and enabling more powerful ambiguity set modeling. For example, moment-based ambiguity sets are more prevalent in DRMDPs compared to RMDPs.

Although RMDPs and DRMDPs differ in certain aspects, their formulations are closely aligned within the scope of this article, and the high-level modeling approaches exhibit substantial parallels. In fact, with appropriate modeling and practical considerations, one formulation might be transformed into the other (e.g., exogenous randomness specified later). Thus, this article mainly investigates the literature from a modeling framework perspective rather than emphasizing the distinctions between RMDPs and DRMDPs. Unless explicitly noted otherwise, RMDPs and DRMDPs will both be referred to as RMDPs in the following for simplicity.

\subsection{Rectangularity}
While the formulation of RMDPs is intuitive and builds on standard MDPs, some properties such as the recursive structure of MDPs do not naturally extend to RMDPs. Specifically, efficiently solving MDPs relies on the Bellman optimality equation, which decomposes the multistage problem into a recursive form. However, extending this to a robust counterpart (i.e., robust Bellman optimality equation) is not trivial, as it accounts for the global worst-case scenario, which can depend on the entire history rather than just the current state. An effective approach to ensuring the computational benefits of RMDPs is to introduce rectangular assumptions. 

Rectangularity is a common assumption in multistage problems and has received great attention from the economics and optimization communities. In economics, the rectangular assumption implies dynamic consistency or time consistency, ensuring the dynamic behavior of a decision-maker is fully determined by his/her preference, rather than conditional preference after each history separately \citep{strotz1973myopia, epstein2003recursive}. In optimization, particularly in stochastic programming (SP), most stochastic dual dynamic programming (SDDP) algorithms require this assumption such that they can safely decompose the multistage problems into $T$ two-stage problems (assuming a planning horizon of $T$), which are then iteratively solved through forward and backward pass \citepalias{philpott2008convergence, zou2019stochastic}. Such decomposition is feasible because the rectangular assumption precludes the existence of an ``inter-temporal budget" for unfavorable outcomes. However, the concept of rectangularity in RMDPs carries additional nuances. First, due to the definition of MDPs, the rectangularity naturally extends to stage-wise and state-action-wise (or state-wise) structures. Second, and more crucially, beyond the independence between inter-temporal and inter-state randomness, the rectangular assumption in RMDPs indicates the independence from external uncertainty (recall that external uncertainty arises from incomplete information about the dynamics). In other words, external uncertainty cannot be changed regardless of the realization in the last stage or state. This latter nuance is also demonstrated in robust dual dynamic programming (RDDP) \citepalias{georghiou2019robust}. While some papers in SDDP and RDDP attempt to work without the rectangular assumption such as restricting the decision rules \citepalias{daryalal2023two, daryalal2024lagrangian}, the rectangular assumption remains important and brings significant convenience for analysis.

As rectangularity is widely applied across various fields and a complete consensus on a common rigorous definition has not yet been reached \citep{xin2021time}, we present several popular definitions to convey the core idea of rectangularity.

\begin{definition}[Rectangularity] 
     For simplicity, we term the collection of uncertainty as an ambiguity set, denoted by $\mathcal{P}$. In terms of mathematical representation, the definition of rectangularity can be categorized into three types:

    (1) \textit{Static form} \citep{nilimRobustControlMarkov2005,le2007robust,wiesemann2013robust}. In this case, the state-wise rectangular ($s$-rectangular) ambiguity set $\mathcal{P}^{\mathcal{S}}$ and state-action-wise rectangular ($(s,a)$-rectangular) ambiguity set $\mathcal{P}^{\mathcal{SA}}$ are time-invariant. Taking the transition probability matrix, for instance, both ambiguity sets can be represented respectively as follows:
    \begin{equation}
        \begin{aligned}
        \mathcal{P}^{\mathcal{S}} &=\left\{ \mathbf{p}\in(\Delta^{|\mathcal{S}|})^{|\mathcal{S}|\times|\mathcal{A}|}\,\middle |\, \mathbf{p}=\bigotimes_{s\in\mathcal{S}}\mathbf{p}_{s},\mathbf{p}_{s}\in\mathcal{P}_{s},\mathcal{P}_{s}\subseteq\Delta^{|\mathcal{S}|\times|\mathcal{A}|}\right\}, \\
        \mathcal{P}^{\mathcal{SA}} &=\left\{ \mathbf{p}\in(\Delta^{|\mathcal{S}|})^{|\mathcal{S}|\times|\mathcal{A}|} \,\middle |\, \mathbf{p}=\bigotimes_{(s,a)\in\mathcal{S}\times\mathcal{A}}\mathbf{p}_{sa},\mathbf{p}_{sa}\in\mathcal{P}_{sa},\mathcal{P}_{sa}\subseteq\Delta^{|\mathcal{S}|}\right\}.
        \end{aligned}
    \end{equation} 
    where $\mathcal{P}_s$ and $\mathcal{P}_{sa}$ are the set of (conditional) transition probability matrices of interest.
    
    (2) \textit{Dynamic form} \citep{iyengarRobustDynamicProgramming2005, epstein2003recursive, iancuTightApproximationsDynamic2015}. This representation is closely related to the notion of time consistency. Recall $\mathbb{P}(\cdot\mid \mathcal{F}_t)$ is the probability measure $\mathbb{P}$ conditioned on the filtration $\mathcal{F}_t$. Let $\mathbb{P}^{+1}(\cdot|\mathcal{F}_{t})$ be the restriction of $\mathbb{P}(\cdot\mid \mathcal{F}_t)$ to $\mathcal{F}_{t+1}$. In this case, the ambiguity set is rectangular if it satisfies:
    \begin{equation}
    	\mathcal{P}_t=\left\{ \int_{\Xi_{t+1}} \mathbb{P}(\cdot\mid \mathcal{F}_{t+1})\mathrm{d}\mathbb{Q} \,\middle |\, \mathbb{P}\in\mathcal{P}, \mathbb{Q}\in\mathcal{P}_{t}^{+1}\right\},
    \end{equation}
    where $\mathcal{P}_{t}^{+1}=\left\{ \mathbb{P}^{+1}(\cdot\mid \mathcal{F}_{t})\mid \mathbb{P}\in\mathcal{P}\right\}$ contains the conditional one-step-ahead measures.

    (3) \textit{Nested form} \citep{shapiro2016rectangular}. Let $f(\xi_{[T]})$ be the cumulative discounted reward for simplicity, where $\xi_{[T]}=(\xi_1,\xi_2,...,\xi_T)$ is the history of the data up to time $T$. In the clear context, we can bridge the equivalence between $\xi$ and transition  $(s, a, s')$. In this case, the rectangular assumption is satisfied if the following equation holds,
    \begin{equation}
    \label{nested form}\sup_{\mathbb{P}\in\mathcal{P}}\mathbb{E}_{\mathbb{P}}\left[f(\xi_{[T]})\right]=\sup_{\mathbb{P}\in\mathcal{P}}\mathbb{E}_{\mathbb{P}}\left[\sup_{\mathbb{P}\in\mathcal{P}}\mathbb{E}_{\mathbb{P}|\xi_{1}}\left[\cdots\sup_{\mathbb{P}\in\mathcal{P}}\mathbb{E}_{\mathbb{P}|\xi_{[T-1]}}[f(\xi_{[T]})]\right]\right],
    \end{equation}
    where $\mathbb{E}_{\mathbb{P}|\xi_{[t]}}[f(\xi_{[T]})]$ is the conditional expectation of $f(\xi_{[T]})$ with respect to $\mathbb{P}\in\mathcal{P}$ given $\xi_{[t]}$.
\end{definition}

The first form is termed ``static" because it lacks a subscript $t$. The static form of rectangularity aims to capture the state-wise or state-action-wise independence with rectangularity, and a more concise representation is given by $\mathcal{P}^{\mathcal{S}} = \bigotimes_{s\in\mathcal{S}}\mathcal{P}_s$ and $\mathcal{P}^{\mathcal{SA}}=\bigotimes_{s\in\mathcal{S},a\in\mathcal{A}}\mathcal{P}_{sa}$. It is noteworthy that the static form is often the most popular representation of rectangularity in RMDP literature. In economics, particularly risk measure community, the dynamic form of the rectangular structure is derived using one-step conditional measures, motivated by the dynamic consistency of the expected utility model. From a scenario tree perspective, this dynamic representation means that we can only consider the leaves that are from specific parent nodes. Recently, \cite{shapiro2016rectangular} proposed a nested-form condition that ensures the rectangularity assumption, directly targeting the decomposability property of the static formulation to establish the equivalence with the dynamic formulation. Obviously, as the supremum operator is taken over time, the right-hand side of \eqref{nested form} is greater than or equal to the left side. With rectangular assumption, the optimization problem of $t$ suffices to take the maximum over the corresponding set of marginal probability measures of the observations $\xi_{[t-1]}$. While there are interesting discussions on the relationship between rectangularity and time-(in)consistency, they are beyond the scope of this review. We refer interested readers to \cite{xin2021time} and the references therein for details.

\begin{remark}
	It is necessary to recognize that these definitions of rectangularity are not simply semantically distinct but differ fundamentally in their underlying interpretations. However, no matter which of the above three definitions, with any one of the three rectangular assumptions, a recursive representation can be derived from the original multistage problem, thereby rendering the formulation tractable. Meanwhile, stage-wise independence stands as a common condition that can satisfy rectangularity as characterized in any formulation (In RMDPs, we always assume the sets of states at different stages are mutually exclusive).
\end{remark}

Returning to the static form of rectangularity, which is most prevalent in RMDPs, we further detail the nature and distinctions between $(s, a)$-rectangular ambiguity sets and $s$-rectangular ambiguity sets. In the context of game theory, these two types of ambiguity sets differ in their assumptions regarding nature's knowledge and commitment \citep{le2007robust}. The $(s, a)$-rectangular ambiguity sets essentially presuppose that nature can observe the decision maker's actions before choosing the worst plausible realization of the transition probabilities. In other words, nature can react to the decision maker's actions and choose the worst-case scenario consistent with the observed actions, often referred to as a \textit{reactive} or \textit{adaptive} nature. In contrast, $s$-rectangular ambiguity sets assume that a weaker nature must commit to a realization of the transition probabilities before observing the decision-maker's actions. Here, nature must choose a realization of the transition probabilities without the benefit of observing the decision maker's actions, ensuring that this realization is the worst case among all possible realizations consistent with the ambiguity sets. While $s$-rectangular and  $(s, a)$-rectangular assumptions are suited to different contexts (e.g., whether can obtain independent transition samples for state-action pair), $s$-rectangular ambiguity sets generally (but not always) result in less conservative solutions yet with a higher computational complexity \citep{wiesemann2013robust, ho2022robust}. The intuition behind this lies in nature's inability to adapt to the decision-maker's actions under $s$-rectangular ambiguity, resulting in a less pessimistic worst-case scenario than that of $(s, a)$-rectangular sets.

The rectangular assumptions have been investigated widely about the tractability of RMDPs, and \cite{Bagnell2001SolvingUM} and \cite{wiesemann2013robust} prove that solving a general (non-rectangular) RMDP is NP-hard. We provide a new proof that is similar to \cite{Bagnell2001SolvingUM} by constructing a non-rectangular RMDP as a conjunctive normal form (CNF) formula yet more straightforward. In addition,  we also show how to add the rectangularity assumptions such that the RMDPs can be solved in polynomial time.

\begin{proposition}[NP-hardness for non-rectangularity]
    \label{np-hard proof}
    Solving the non-rectangular RMDPs is NP-hard.
\end{proposition}

\proof{of Proposition \ref{np-hard proof}}
We prove the NP-hardness of solving non-rectangular RMDPs by reduction from the Conjunctive Normal Form Satisfiability (CNF-SAT) problem. A Boolean formula is in CNF if it is a conjunction $\wedge$ of clauses (or a single clause), where each clause is a disjunction $\vee$ of literals (or a single literal). The CNF-SAT is the decision problem of determining if there exists an assignment of truth values (``True" or ``False") to the variables such that the given formula evaluates to ``True". For instance, the formula $(x_1\vee \neg x_2)\wedge(\neg x_1 \vee x_3)$ is in CNF, with two clauses $C_1=x_1\vee \neg x_2,C_2=\neg x_1\vee  x_3$, and three literals (variables) $x_1,x_2,x_3$, where $\neg$ denotes the negation. If we assign $x_1,x_2,x_3$ all as ``True", the formula evaluates to ``True", thereby providing a solution.  While the problem may seem straightforward, finding a solution is NP-hard. Next, we will demonstrate how to bridge the connection between CNF-SAT and non-rectangular RMDPs.

Let a CNF formula be $\phi=C_1\wedge C_2\wedge \cdots\ C_m$ over variables $x_1,x_2,...x_n$. Given the CNF formula, we construct a non-rectangular RMDP $\mathcal{M}=<\mathcal{S},\mathcal{A},\mathbf{p}>$ without discounting and reward randomness as follows:

\begin{enumerate}
	\item The state space $\mathcal{S}$ consists of $n+m$ states: For each clause $C_i$, we create a state $s_i$. For each variable $x_j$ in clause $C_i$, we create state $s_j$ which is reachable from state $s_i$.
	\item The action space $\mathcal{A}$ consists of the single action of assigning truth values to variables in clauses.
	\item The transition probability $p_{ij}$ represents the probability of assigning ``True" to variable $x_j$ in clause $C_j$. Consequently, $1-p_{ij}$ is the probability of assigning ``False". Note that the action corresponds to the assignment, thus $p_{ij}$ implicitly depends on the action, consistent with the notion of transition probability $p_{sas'}$ in RMDPs.
\end{enumerate}

Consider a deterministic policy $\pi:\mathcal{S}\rightarrow \mathcal{A}$. Under $\pi$, each $p_{ij}$ becomes binary, either 0 or 1: if ``True" is assigned to variable $x_j$ in clause $C_i$, the state transitions from $s_i$ to $s_j$ with probability 1. The constructed RMDP $\mathcal{M}$ is non-rectangular because the transition probabilities $p_{ij}$ are not independent across states. Specifically, an assignment of a truth value to $x_j$ in one clause (state) affects all clauses containing $x_j$ simultaneously. This creates a coupled uncertainty structure where the choice of action in one state constrains the possible transitions in other states.

Suppose there exists an algorithm $\mathbf{A}$ that solves $\mathcal{M}$. We can use $\mathbf{A}$ to solve the original CNF-SAT problem as follows:
\begin{enumerate}
	\item Construct $\mathcal{M}$ from the given CNF formula $\phi$ as described above.
	\item Apply algorithm $\mathbf{A}$ to find an optimal policy $\pi^*$ for $\mathcal{M}$.
	\item Interpret $\pi^*$ as a satisfied assignment for $\phi$: if $\pi^*(s_i)$ transitions to $s_j$, set $x_j$ to ``True" in $C_i$; otherwise set in to ``False".
\end{enumerate}
If algorithm $\mathbf{A}$ operates in polynomial time, this procedure would enable CNF-SAT to be solved in polynomial time. However, CNF-SAT is a known NP-complete problem. Therefore, since the reduction is polynomial-time and a solution to the RMDP provides a solution to CNF-SAT, it follows that solving non-rectangular RMDPs is NP-hard. 

Now, we demonstrate how to leverage the rectangular assumption that makes constructed $\mathcal{M}$ solvable in polynomial time. In rectangular RMDPs, the assumption that each state is visited only once, and each visit is to a unique state is required. In the context of CNF-SAT reduction, if this assumption and stage-wise (clause-wise) independence are satisfied, the global consistency constraint is essentially relaxed. Formally, we remove the constraint that $x_j^i=x^k_j$ for all clauses $i$ and $k$ containing variable $j$. This relaxation makes each state (clause) independent, creating a rectangular RMDP where each state can be optimized separately. However, the solution may not be valid for the original non-rectangular RMDP or CNF-SAT problem due to potential inconsistencies. This availability of such decomposition highlights the significance of the rectangular assumption in RMDPs.
\endproof

Finally, we briefly present several novel rectangular concepts to close this section. \citetalias{mannor2016robust} propose the concept of $k$-rectangularity, which restricts the cardinality of deviations to no more than $k$ to maintain tractability. Additionally, \citetalias{gohDataUncertaintyMarkov2018} and \cite{goyalRobustMarkovDecision2023} represent transition probabilities as linear combinations of $r$ factors and introduce a new assumption termed $r$-rectangularity. These novel rectangular concepts render RMDPs tractable without relying on traditional rectangular assumptions (e.g., state-wise independence), at the expense of representation flexibility or independence of factors. More details about these novel rectangular concepts will be specified in Section \ref{sec6}. Whether using traditional rectangular assumptions or exploring novel frameworks, the principle of ``no free lunch'' applies. Ensuring the tractability of RMDPs inherently requires the adoption of certain (potentially strong) assumptions.

\section{Parametric RMDPs}\label{sec3}

Parametric RMDPs are among the earliest studied RMDPs \citep{satia1973markovian}, which assume that transition probabilities are controlled by known parameters. Through predetermining supports of parameters or distribution families, parametric RMDPs transfer the uncertainty sets of transition probabilities into uncertainty sets of the parameters. Although the study of parametric RMDPs has undergone significant changes, they are grouped under the same category in this paper from a high-level perspective. These major shifts in research have led to the exploration of new approaches and techniques within the field of parametric RMDPs.

As a pioneering work, \cite{satia1973markovian} propose an unprecedented formulation, known as MDPs with imprecisely known parameters (MDPIPs), where uncertainty sets comprise individually bounded transition probabilities with given lower and upper bounds. The uncertainty sets $\mathcal{P}$ are defined by:
\begin{equation}
    \label{MDPIPs sets}
    \mathcal{P}\triangleq\left\{\mathbf{p}\in (\Delta^{|\mathcal{S}|})^{|\mathcal{S}|\times|\mathcal{A}|}\left|
    \begin{array}{l}
        p^{L}_{sas'}\leq p_{sas^{\prime}}\leq p^{U}_{sas'},\\
        \textstyle\sum_{s'}p_{sas^{\prime}} = 1,\\
        p_{sas^{\prime}}\geq 0,~\forall s, s'\in \mathcal{S}, a\in \mathcal{A}
    \end{array}
    \right\}\right.,\\
\end{equation}
where $p^{L}_{sas'}$ and $p^{U}_{sas'}$ are the predetermined lower and upper bounds, respectively. The uncertainty set $\mathcal{P}$ in \eqref{MDPIPs sets} builds a polyhedron. Particularly, we argue that the authors implicitly assume the $(s, a)$-rectangular assumption. They show that there exists a pure, stationary policy that is optimal from a game theory perspective because nature always can observe the decision-maker's decision before making its decision under $(s, a)$-rectangularity. Consequently, the authors develop a modification of the policy iteration (PI) algorithm described by \cite{howard1960dynamic}, where the objective of the policy evaluation procedure is minimizing the expected return over any policy, instead of maximizing in the original PI algorithm. They prove the proposed algorithm can find an $\epsilon$-optimal policy in a finite number of iterations. Subsequently, \cite{white1994markov} study the same formulation as \cite{satia1973markovian}, but propose a new PI algorithm that leverages the LP techniques to accelerate the convergence effectively.

Following the same line of research, \cite{givan2000bounded} investigate a special subclass of MDPIPs, termed as \textit{bounded-parameter} MDPs (BMDPs). A BMDP, denoted as $\mathcal{M}_{\updownarrow}=<T,\mathcal{S},\mathcal{A},R_{\updownarrow},\mathbf{p}_{\updownarrow},\gamma>$, is similar with the general MDP, but $\mathbf{p}_{\updownarrow}$ and $R_{\updownarrow}$ are closed real intervals instead of fixed real values. For instance, given a transition tuple $(s,a,s')$, the probability $p_{sas'}$ in MDPs is a real value within range $[0,1]$. However, in BMDPs,  $p_{\updownarrow,sas'}$ is expressed as a closed interval of the form $[p_{\downarrow,sas'},p_{\uparrow,sas'}]$ where $0\le p_{\downarrow,sas'}\le p_{\uparrow,sas'} \le 1$ are known parameters, similar to $p^{L}_{sas'}, p^{U}_{sas'}$ in \eqref{MDPIPs sets}. Thus, $\mathcal{M}_{\updownarrow}$ essentially represents a set of exact MDPs, where an MDP $M\in \mathcal{M}_{\updownarrow}$ if $M=<T,\mathcal{S},\mathcal{A}, R',\mathbf{p}',\gamma>$ and $R(s, a,s')\in R_{\updownarrow,sas'}, p _{sas'}\in p_{\updownarrow,sas'}$ for all $s,s'\in \mathcal{S}, a\in \mathcal{A}$.

The construction of $\mathcal{M}_{\updownarrow}$ highlights the distinction with other parametric RMDPs, as the key elements in BMDPs are different exact MDPs, rather than uncertain transition probabilities. Given a well-defined $\mathcal{M}_{\updownarrow}$, BMDPs focus on an interval value function $V_{\updownarrow}$ which is a mapping from states to closed real intervals. For any state $s\in\mathcal{S}$, $V_{\updownarrow}(s)=[V_{\downarrow}(s),V_{\uparrow}(s)]$, where $V_{\downarrow}(s)$ and $V_{\uparrow}(s)$ are some real values. Although the element number of $\mathcal{M}_{\updownarrow}$ is generally infinite, the authors show that the attention can be restricted within a particular MDPs family as the \textit{order-maximizing MDPs}. For instance, given an arbitrary state permutation $\mathcal{O}$, its transition probabilities are assigned as the following threshold structure given $s,a$:
\begin{equation}
\label{threshold_BMDP}
	p_{sas_i}=\begin{cases}
p_{\uparrow,sas_i}, \text{ if }i\le r-1\\
p_{\downarrow,sas_i}, \text{ if }i\ge r\\
\end{cases}
\end{equation}
Where $s_i$ denotes the $i$-th state in the permutation $\mathcal{O}$ and $r=\arg\max_j \sum_{i=1}^{j-1}p_{\uparrow,sas_i}+\sum_{i=j}^{|\mathcal{S}|}p_{\downarrow,sas_i}$. In other words, an order-maximizing MDP aims to assign more transition probabilities into early states in permutation $\mathcal{O}$. Building upon this class of MDPs, the authors develop interval value iteration where $V_{\downarrow}$ is obtained in terms of a permutation with increasing value and $V_{\uparrow}$ is computed by a decreasing order. By restricting the class of MDPs to the order-maximizing ones, the computational cost is significantly reduced compared to the traditional MDPIPs with a polyhedron uncertainty set.

More recently, \citetalias{delimpaltadakis2023interval} also leverage the threshold structure like \eqref{threshold_BMDP} to expand the scope of BMDPs into continuous action spaces. Distinct from simple real-valued intervals, the authors introduce two functions $f^L, f^U:\mathcal{S}\times\mathcal{A}\times\mathcal{S}\to[0,1]$ to represent the lower and upper bounds on the transition probability with respect to transition $(s, a, s^{\prime})$. Due to strong duality and the ordering threshold structure, they successfully convert the original max-min formulation into a maximization problem. Furthermore, they develop an efficient value iteration algorithm based on this reformulation when $\mathcal{A}$ is a polytope.

Besides through the lens of specific ranges of individual transition probability, \cite{osogami2012robustness} considers parametric RMDPs with $s$-rectangular uncertainty sets specified by a factor, $0<\alpha<1$, which determines the possibly maximal value $\frac{1}{\alpha}\hat{p}_{sas^{\prime}}$ of each transition probability $p_{sas'}$ based on the nominal values $\hat{p}_{sas^{\prime}}$: 
\begin{equation}
    \mathcal{P}\triangleq\left\{\mathbf{p}\in (\Delta^{|\mathcal{S}|})^{|\mathcal{S}|\times|\mathcal{A}|}\,\left|\, 0\leq p_{sas^{\prime}}\leq\frac{1}{\alpha}\hat{p}_{sas^{\prime}}\text{ and }\textstyle\sum_{s'\in\mathcal{S}}p_{sas^{\prime}}=1\ \forall s, s'\in\mathcal{S}, a\in \mathcal{A} \right\}\right. .
\end{equation}
Instead of explicit intervals for each probability, this approach requires only a single parameter $\alpha$, which makes it more convenient to construct uncertainty sets and achieve better generalization performance when the decision-maker has highly limited information about transitions. Furthermore, one can interpret $\alpha$ as the level of robustness or risk-aware confidence level. From this perspective, the author establishes the equivalence between RMDPs and risk-sensitive MDPs, where the risk confidence level is set to be $1-\alpha$. This relationship provides further insights to understand and construct uncertainty sets.

Different from the purely polytopic uncertainty sets, \cite{wiesemann2013robust} consider uncertainty sets $\mathcal{P}$ where transition probabilities follow an affine function of a vector $\theta$:
\begin{equation}\label{eq:affine_set}
    \mathcal{P}\triangleq\left.\left\{\mathbf{p}\in (\Delta^{|\mathcal{S}|})^{|\mathcal{S}|\times|\mathcal{A}|}\,\right|\,\exists\ \theta\in \Theta\text{, such that }\mathbf{p}_{sa}:=P^\theta(\cdot\mid s,a)\ \forall(s,a)\in \mathcal{S}\times \mathcal{A}\right\},
\end{equation}
where $\Theta$ is a subset of $\mathbb{R}^q$ and $P^\theta(\cdot\mid s,a),~\forall(s,a)\in \mathcal{S}\times \mathcal{A}$ is an affine function from $\Theta$ to $[0,1]$ that satisfies $P^\theta(\cdot\mid s,a):=k_{sa}+K_{sa}\theta$ for some $k_{sa}\in \mathbb{R}^{|\mathcal{S}|}$ and $K_{sa}\in \mathbb{R}^{|\mathcal{S}|\times k}$. The affine function form implies that the conditional transition probability is the outcome of a linear regression, where the coefficients are $(K_{sa},k_{sa})$ regarding the parameters $\theta$. Note that $\theta$ controls the ambiguity sets essentially, which allows to condense all ambiguous information and parameters in the set $\Theta$, e.g.,
\begin{equation}
    \label{wieseman_support}
    \Theta\triangleq\left\{\theta\in\mathbb{R}^q\colon\theta^\top O_l\theta+o_l^\top\theta+\omega\geq0,~\forall l=1,\ldots,L\right\},
\end{equation}
where $O_l$ is a $k\times k$ matrix, $O_{l}\preceq0$, $\omega$ is the constant and $L$ is the number of constraints applied to $\theta$. We assume that $\Theta$ has a nonempty interior, that is, none of the parameters in $\Theta$ is fully explained by the others. The benefit of this construction defined in \eqref{wieseman_support} is that the quadratic inequalities in $\Theta$ with negative semidefinite matrices can be reformulated as second-order cone constraints, allowing RMDPs to be expressed as second-order cone programs (SOCPs) or semidefinite programs (SDPs), which are convex problems and can be solved efficiently with interior-point methods. The authors assume that $\Theta$ is bounded and has a nonempty interior, which means none of the parameters in $\Theta$ is fully explained by the others. With the $s$-rectangular assumption, they show that the robust Bellman Optimality equation holds with the optimization procedure proceeds for each state $s\in \mathcal{S}$ and such operators are contraction mappings, ensuring the unique optimal value-to-go function. Furthermore, the authors provide the time complexity of policy evaluation and policy improvement routines, respectively.

\citetalias{black2022robust} apply RMDPs for newsvendor problems, where the demand is exogenous and state transitions are essentially driven by demand realizations (underlying dynamics). Thus, the transition probability matrix is not time-invariant, and given the inventory position, the conditional transition probabilities can be equivalently derived by the distribution of demands. In this case, the authors investigate RMDPs where the ambiguity sets are limited to distributions within the same parametric family, focusing on identifying the worst-case parameters rather than the entire distribution. This approach aims to alleviate challenges with moment estimation in non-parametric RMDPs, as well as the issue of overly conservative solutions of non-parametric methods that can arise when parametric distributions are fitted well enough. Let $f_{\xi_{sa}}$ be the probability mass function of exogenous random variable $\xi_{sa}$ which is controlled by the parameters $\boldsymbol{\theta}_{sa}=(\theta_{sa_1},\cdots,\theta_{sa_k})\in \mathbb{R}^k$. The authors assume that the next state $s_{t+1}$ is specified by a simple, known function $g$ of $\xi_{sa}$ as $s_{t+1}=g(\xi_{sa}|s, a)$. For example, in inventory management, states often denote the initial inventory positions and actions are the order quantities. Consequently, one possible function is $g=s+a-\hat{\xi}_{sa}$ with the demand realization $\hat{\xi}_{sa}$. As we discussed in REMARK \ref{probability_remark}, for each $(s, a)\in \mathcal{S}\times\mathcal{A}$, conditional transition (discrete) distribution $\mathbf{p}_{sa}$ can be uniquely represented as the distribution of $\xi_{sa}\in \Xi_{sa}$. More precisely, each transition probability $p_{sas^{\prime}}$ for the transition tuple $(s,a,s^{\prime})$ can be computed as
\begin{equation}
\label{parametric_structure}
    p_{sas^{\prime}} = P^{\boldsymbol{\theta}_{sa}} (s^{\prime}\,|\,s,a)=P^{\boldsymbol{\theta}_{sa}}(g(\xi_{sa}\,|\,s,a)=s^{\prime})=\sum_{\xi\in \Xi_{sa}(s^{\prime})} f_{\xi_{sa}}(\xi\,|\,\boldsymbol{\theta}_{sa}) ,
\end{equation}
where $\Xi_{sa}(s^{\prime})$ is the support of $\xi$ with $(s,a,s^{\prime})$ transition. Because the transition probability matrix $\mathbf{p}$ is uniquely specified by $\boldsymbol{\theta}$, the authors simplify ambiguity sets for $\boldsymbol{\theta}$, rather than $\mathbf{p}$, to maintain the equivalence structure as in \eqref{parametric_structure}. They consider ambiguity sets with $s$-rectangular structure $\Theta=\bigotimes_{s\in \mathcal{S}}\Theta_s$ where $\Theta_s\subseteq\mathbb{R}^{k},\mathrm{~}\forall s\in\mathcal{S}$,
and then reformulate the RMDP as:
\begin{equation}
    \max_\pi\min_{\boldsymbol{\theta}\in\Theta}\mathbb{E}_{\mathbf{p}}^\pi\left[\left.\sum_{t=0}^\infty\gamma^tR(s_t,a_t,s_{t+1}) \,\right| \,\boldsymbol{\theta}\right].
\end{equation}

In this paper, the authors also demonstrate how to construct the uncertainty set in terms of $\boldsymbol{\theta}$ when finite samples from the true distribution of $\xi$ are accessible. By standard results in maximum likelihood theory \citep{millar2011maximum}, given the log-likelihood function $\ell(\cdot)$ for observed data, a Maximum Likelihood Estimation (MLE) $\hat{\boldsymbol{\theta}}_{sa}$ of true parameter $\boldsymbol{\theta}^0_{sa}$ satisfies:
\begin{equation*}
	\sum_{a\in\mathcal{A}}\left(\hat{\boldsymbol{\theta}}_{sa}-\boldsymbol{\theta}_{sa}^0\right)^TI_{\mathbb{E}}\left(\boldsymbol{\theta}_{sa}^0\right)\left(\hat{\boldsymbol{\theta}}_{sa}-\boldsymbol{\theta}_{sa}^0\right)\thicksim\chi_{k|\mathcal A|}^2,
\end{equation*}
where $I_{\mathbb{E}}\left(\boldsymbol{\theta}^0_{sa}\right)=\left(-\mathbb{E}_{\xi_{sa}}\left[\frac\partial{\partial\theta^0_{sai}\partial\theta^0_{saj}}\ell({\boldsymbol{\theta}}^0_{sa})\right]\right)_{i,j=1,...,k}$ is the expected Fisher information matrix and $\chi_{k|\mathcal A|}^2$ is the chi-squared distribution with $k|\mathcal A|$ degrees of freedom. As the  asymptotic equivalence holds,  the authors use $\hat{\boldsymbol{\theta}}_{sa}$ estimated via historical data to approximate $\boldsymbol{\theta}^0_{sa}$ and construct a $1-\alpha$ confidence level set given by:
\begin{equation}
    \Theta_{s}=\left\{\boldsymbol{\theta}_s\in\mathbb{R}^{|\mathcal{A}|}\times\mathbb{R}^k\,\middle|\,\sum_{a\in\mathcal{A}}\left(\hat{\boldsymbol{\theta}}_{sa}-\boldsymbol{\theta}_{sa}\right)^TI_{\mathbb{E}}\left(\hat{\boldsymbol{\theta}}_{sa}\right)\left(\hat{\boldsymbol{\theta}}_{sa}-\boldsymbol{\theta}_{sa}\right)\leq\chi_{k|\mathcal{A}|,1-\alpha}^2\right\}\quad \forall s\in \mathcal S.
\end{equation} The uncertainty set $\Theta_{s}$ contains all parameter vectors $\boldsymbol{\theta}_s$ that are ``close enough" to the MLE estimate, where ``close enough" is defined such that the sum of the distances (i.e., the quadratic term weighted by the Fisher Information matrix) over all actions must be less than or equal to the $(1-\alpha)$ quantile of $\chi_{k|\mathcal A|}^2$. Considering the non-linearity of $\mathbf{p}$ as functions of the parameters, they replace the uncertainty set $\Theta_{s}$ with a discretized set $\Theta_{s}^{\prime}$ and reformulate the Bellman update \eqref{bellmaneuqation} as a linear program with $|\Theta_{s}^{\prime}|+1$ constraints.

\section{Moment-Based RMDPs}\label{sec4}

The motivations behind parametric RMDPs are relatively intuitive, however, determining reasonable ranges for parameter values in practice remains a significant challenge. Moreover, most parametric approaches intend to impose constraints on each transition probability $p_{sas'}$ such that the uncertainty sets are excessively large. While this ensures general applicability, it often leads to overly conservative or impractical solutions. On the other hand, it is more practical to assume that decision-makers have confidence in some probabilistic information regarding the dynamics. Under such circumstances, moments of the distribution, such as the mean and variance, are typically more accessible and reliable than precise parameter ranges, as they can be readily inferred from historical data.

Since the seminal paper by \cite{scarf1957min} marks the inception, moment-based ambiguity sets have been extensively studied in the literature of RO and DRO for decades \citepalias{delage2010distributionally, wiesemann2014distributionally, bertsimas2019adaptive, ninh2021robust}. While this concept has also been extended to RMDPs, the literature on RMDPs employing moment-based ambiguity sets is considerably limited compared to RO and DRO. This relative scarcity can be attributed to the inherent nature of RMDPs, where the dynamics are typically represented as environmental processes that drive state transitions, rather than as explicit random variables with explicitly interpretable moment information. For instance, in standard RMDPs, the transition probability matrix is defined on a probability simplex. Given state-action pair $(s, a)$, it is challenging to interpret the mean or variance of a vector (i.e., a discrete conditional distribution).

Nevertheless, moment-based ambiguity sets can sometimes be tailored well to RMDPs where the randomness of the decision problem is exogenous, namely, unrelated to states and actions. It is an interesting concept discussed as \textit{exogenous process} recently in model-based RL algorithms \citep{madeka2022deep, pmlr-v202-sinclair23a}. Besides inventory management discussed above, numerous applications exhibit characteristics of exogenous processes, including appointment scheduling, airline revenue management, and resource allocation. While it may not be natural to understand and estimate the moments of transition probabilities in general RMDPs or DRMDPs, it is important to carefully include the moment-based information in the construction of ambiguity sets in many applications.

\cite{yangDynamicGameApproach2018} study the RMDPs where the transition is controlled by state $s$, action $a$, stochastic disturbance $\boldsymbol{\xi}$, and a measurable function $f$ \footnote{More precisely, this formulation is more related to the concept of optimal control (OC). However, under mild assumptions and appropriate settings, we can reformulate an OC as an MDP, regarding significant overlap in the key techniques utilized by both methodologies.}:
\begin{equation*}
	s_{t+1}=f(s_t, a_t,\boldsymbol{\xi}_t) \quad \forall t\in[T].
\end{equation*}
In this case, the author supposes that estimates of the mean, $\boldsymbol{m}_t\in \mathbb{R}^{q}$, and covariance matrix, $\boldsymbol{\Sigma}_t\in \mathbb{R}^{q\times q}$, of the disturbance $\boldsymbol{\xi}_t\sim \mu_t$ are the only available information for each stage $t$. With the stage-wise independent assumption implicitly, the moment-based ambiguity sets are modeled over time as 
\begin{equation}
    \label{momentambiguity}
    \mathcal{P}_t\triangleq\left\{\mu_t\in\mathfrak{M}(\Xi_t,\mathcal{F}^t)~\left|~\begin{array}{l}
    \mu_t(\Xi_t)=1\\
    |\mathbb{E}_{\mu_t}[\boldsymbol{\xi}_t]-\boldsymbol{m}_t|\leq \boldsymbol{\theta}_t\\
    \mathbb{E}_{\mu_t}[(\boldsymbol{\xi}_t-\boldsymbol{m}_t)(\boldsymbol{\xi}_t-\boldsymbol{m}_t)^{\top}]\preceq \beta_t\boldsymbol{\Sigma}_t
    \end{array}\right\}\right.,
\end{equation}
where $\boldsymbol{\theta}_t\in\mathbb{R}^k$ and $\beta_t\geq 1$ are given constants that depend on the confidence in the estimates $\boldsymbol{m}_t$ and $\boldsymbol{\Sigma}_t$. From a geometric perspective, the three constraints in \eqref{momentambiguity} denote the following: (\textit{i}) the support of $\boldsymbol{\xi}_t$ is $\Xi_t$; (\textit{ii}) the mean of $\boldsymbol{\xi}_t$ lies in a ball of size $\boldsymbol{\theta}_t$; and (\textit{iii}) the centered second-moment matrix of $\boldsymbol{\xi}_t$ lies in a positive semidefinite cone. Thus, this ambiguity models how likely $\boldsymbol{\xi}_t$ is to be close to the estimate $\boldsymbol{m}_t$ in terms of the weighted correlation matrix estimate ${\beta_t} \boldsymbol{\Sigma}_t$. The parameters $\boldsymbol{\theta}_t$ and ${\beta_t}$ allow for adjusting the size of the ambiguity set based on the confidence in the estimates.

Notice that solving an RMDP with an ambiguity set (\ref{momentambiguity}) involves an infinite-dimensional maximin optimization problem, which is generally computationally intractable. To overcome the challenge, the author proposes a dual reformulation of the inner minimization problem, which is a semi-infinite maximal optimization program. Leveraging the results from \cite{lasserre2009moments}, it is shown that no duality gap exists in this reformulation. Consequently, by substituting the inner problem in the Bellman equation without sacrificing optimality, the author develops the dual Bellman equation which can be shown to be concave with some mild assumptions, e.g., measurable function $f$ is affine, and state space $\Xi_t$ is convex and compact. 

\citetalias{song2023decision} consider a decision-dependent epidemic control problem where the transition probabilities depend not only on the stochastic epidemiological processes but also on control manners implemented by the policy-maker. Thus, the ambiguity sets are endogenous (or decision-dependent), i.e., the transition probabilities depend on the action $a$ and state $s$. With $(s, a)$-rectangular assumption implicitly, the authors construct ambiguity sets for each $(s, a)$ pair:
\begin{equation}
    \label{decisiondependent}
    \mathcal{P}_{sa}\triangleq\left\{\mu_{sa}\in \mathfrak{M}(\mathcal{S},\mathcal{F}^{\mathcal{S}})\mid\mathbf{p}_{sa}\sim \mu_{sa},~\boldsymbol{\theta}_{sa}^{L}\leq\mathbb{E}_{\mu_{sa}}[\mathbf{p}_{sa}]\leq \boldsymbol{\theta}_{sa}^{U}\right\},
\end{equation}
where $\boldsymbol{\theta}_{sa}^{L}$ and $\boldsymbol{\theta}_{sa}^{U}$ are the lower bound and upper bound vectors of the mean of transition probabilities, and $\mathcal{F}^{\mathcal{S}}$ is the corresponding $\sigma$-algebra when sample space is $\mathcal{S}$. Facing the same issue that RMDPs with moment-based ambiguity set are generally intractable, \citetalias{song2023decision} first relax the hard constraint $\boldsymbol{\theta}_{sa}^{L}\leq\mathbb{E}_{\mu_{sa}}[\mathbf{p}_{sa}]\leq \boldsymbol{\theta}_{sa}^{U}$ into soft constraints:
\begin{equation*}
	\int \mathbf{p}_{sa} \mathrm{d}\mu_{sa}(\mathbf{p}_{sa})-\boldsymbol{\theta}_{sa}^{U}\le \boldsymbol{x},\quad \boldsymbol{\theta}_{sa}^{L}-\int\mathbf{p}_{sa} \mathrm{d}\mu_{sa}(\mathbf{p}_{sa})\le \boldsymbol{x},
\end{equation*}
and adjust the objective function (i.e., the Bellman equation) by penalizing constraint violations with $k \mathbf{1}^T\boldsymbol{x}$, where $k$ represents a user-specified penalty coefficient. Subsequently, the authors apply the standard Lagrangian dualization approach to obtain the dual Bellman equation. Slightly different from direct dualization, penalty coefficient $k$ offers greater flexibility and allows the incorporation of expert knowledge into the model. Furthermore, the authors consider $\boldsymbol{\theta}_{sa}^{L}$ and $\boldsymbol{\theta}_{sa}^{U}$ can be expressed by the linear functions of action $a$, e.g.,
\begin{equation*}
	\begin{aligned}
		\theta_{sa}^U(s')=\rho_0(s') + \sum_{i=1}^{N_a}\rho_i^s(s')a_i \quad \forall s'\in \mathcal{S},\\
		\theta_{sa}^L(s')=\eta_0(s') + \sum_{i=1}^{N_a}\eta_i^s(s')a_i \quad \forall s'\in \mathcal{S},
	\end{aligned} 
\end{equation*}
where $a_i$ denotes the $i$-th dimension of action $a$ and $N_a$ is the total dimension. Coefficients $\boldsymbol{\rho}^s$ and $\boldsymbol{\eta}^s$ can be obtained by linear regression or machine learning approaches, and $\rho_i^s(s')$ (or $\eta_{i}^s(s')$) is the coefficient of transition $(s, a,s')$ where $i$-th dimension of action $a$ is $a_i$. Although the linear decision rule overcomes the challenge of dimensionality, it introduces the bilinear terms in the dual Bellman equation (i.e., $\boldsymbol{\theta}_{sa}^U$ or $\boldsymbol{\theta}_{sa}^L$ with the Lagrangian multipliers). The authors adopt McCormick envelope relaxation \citep{mccormick1976computability} and exact unary expansion \citep{gupte2013solving}, leading to a mixed integer programming formulation to represent each Bellman equation.

Rather than purely utilizing the moment information, \citetalias{yu2015distributionally} and \citetalias{chen2019distributionally} both investigate the generalized-moment-based ambiguity sets, which are also called \textit{lifted} ambiguity sets \citep{wiesemann2014distributionally}. Generalized-moment-based ambiguity sets involve probability constraints on the support and expectation constraints on the generalized moment of the ambiguous distributions. Before leveraging the powerful modeling capabilities of lifted ambiguity sets, two critical conditions must be met. Firstly, the uncertain parameters across different states must be independent, namely, $s$-rectangular. Secondly, the admissible state-wise ambiguity set $\mathcal{P}_s$ for each $s\in\mathcal{S}$ must be representable as the union of marginal distributions across all joint distributions of $(\mathbf{p}_s,\tilde{n}_s)$. Here $\tilde{n}_s\in \mathcal{N}_{s}$ is a one-dimensional auxiliary random variable that denotes a scenario. In essence, we could view $\tilde{\mathcal{S}}=\mathcal{S} \times\mathcal{N}_s$ as the new state space in which the original state space incorporates the scenario information.
 
In \citetalias{yu2015distributionally}, the RMDPs with lifted ambiguity sets for each state $s\in\mathcal{S}$ can be expressed as follows:
\begin{equation}
    \label{gm}
    \begin{aligned} 
    \mathcal{P}_s \triangleq 
    \left\{\mu\in \mathfrak{M}(\tilde{\mathcal{S}}, \mathcal{F}^s)\left|
    \begin{array}{lll} 
    &\left(\mathbf{p}_s, \tilde{n}_s\right) \sim \mu \\ \vspace{1ex}
    &\mathbb{E}_{\mu}\left[\mathbf{p}_s \mid \tilde{n}_s \in \mathcal{N}_j\right]=\theta_j & \forall j \in\left[J_s\right] \\ \vspace{1ex}
    &\mu\left[\mathbf{p}_s \in \mathcal{D}_n \mid \tilde{n}_s=n \right]=1 & \forall n \in\mathcal{N}_{s} \\
    &\mu\left[\tilde{n}_s=n\right]=\omega_n\ \text{for some } \boldsymbol{\omega} \in \Delta^{|\mathcal{N}_{s}|} & \forall n \in\mathcal{N}_{s}\\
    &\theta_j \in \mathcal{U}_j & \forall j \in \left[J_s\right]
    \end{array}
    \right\}\right.
    \end{aligned},
\end{equation}
where for each $j\in[J_{s}]$, the set $\mathcal{N}_j\subseteq\mathcal{N}_s$ is a subset of scenarios, and sets $\mathcal{D}_n\subseteq \Delta^{|\mathcal{S}|}$. The notion of $J_s$ is the cardinality of possible scenarios with state $s$. The constraints in \eqref{gm} denote the following: ($i$) conditional mean of $\mathbf{p}_s$ is $\theta_j$ in terms of $\tilde{n}_s \in \mathcal{N}_j$; ($ii$) support constraints $\mathbf{p}_s\in \mathcal{D}_n$ for each scenario $\tilde{n}_s=n$; ($iii$) probability mass (or weight) for each scenarios, where $\boldsymbol{\omega}$ is a probability simplex; and ($iv$) additional constraints on parameters.

In \cite{chen2019distributionally}, the authors propose the RMDPs with lifted ambiguity sets that are a hybridization of a generalized-moment-based ambiguity set and a statistical metric ambiguity set, by adding discrepancy constraints into $(\ref{gm})$:
\begin{equation}\label{discrepancy constraints}
	\mathbb{E}_{\mu}\left.\left[g_{j \tilde{n}_s}(\mathbf{p}_s)~ \right| ~\tilde{n}_s \in \mathcal{N}_j\right] \leq \beta_j  \quad\forall j \in\left[J_s\right],
\end{equation}
and replacing parameter constraints $\{\theta_j \in \mathcal{U}_j\ \forall j \in \left[J_s\right]\}$ with 
\begin{equation}
	(\theta_j,\beta_j)\in \mathcal{U}_j  \quad\forall j \in\left[J_s\right].
\end{equation}
The function $g_{j \tilde{n}_s}$ in \eqref{discrepancy constraints} is convex lower semi-continuous, and the authors explore and formulate $g_{j \tilde{n}_s}$ as the Wasserstein distance. Additionally, they show that this enhanced formulation can be effectively modeled by algebraic modeling packages \citepalias{chen2020robust} and solved by off-the-shelf commercial solvers.

\section{Discrepancy-Based RMDPs}\label{sec5}

In addition to estimating (partial) moment information, it is possible for the decision-maker to obtain a nominal/reference distribution to approximate the underlying probability distribution, either directly from historical data or domain knowledge. If the available data or domain knowledge is reasonably reliable, such as high-quality data and recognizable patterns of uncertainty, it is reasonable to believe that the discrepancy between the nominal distribution and the true distribution is sufficiently small. Such a discrepancy-based ambiguity set can be represented by the following generic form:
\begin{equation}
    \label{discrepancy_ambiguity}
    \mathcal{P}\triangleq\left\{\mathbb{P}\in\mathfrak{M}\left(\Xi,\mathcal{F}\right)\left|~D(\hat{\mathbb{P}},\mathbb{P})\leq\theta\right\}\right.,
\end{equation}
where $\hat{\mathbb{P}}$ denotes the nominal probability measure, and $D(\cdot,\cdot):\mathfrak{M}\left(\Xi,\mathcal{F}\right)\times\mathfrak{M}\left(\Xi,\mathcal{F}\right)\rightarrow \mathbb{R}_+\cup \{\infty\}$ is a function that measures the discrepancy between two probability measures $\mathbb{P},~\hat{\mathbb{P}}\in \mathfrak{M}\left(\Xi,\mathcal{F}\right)$. The parameter $\theta\in[0,\infty]$, also called the level of robustness,  limits the maximum discrepancy, thereby controlling the size of the ambiguity set. As before, we sightly abuse the notation $\mathcal{P}$ to denote the ambiguity set of corresponding $\mathbf{p}$ or $\mu$ with respect to the probability measure $\mathbb{P}$ for simplicity. In line with \cite{rahimian2019distributionally}, we refer to RMDPs with ambiguity sets like (\ref{discrepancy_ambiguity}) as discrepancy-based RMDPs, and focus on several prevalent ambiguity sets in this section.

\subsection{Norm Ambiguity}

Norm ambiguity sets, typically constructed on a probability simplex, measure the discrepancy between two distributions regarding various norms. While this construction appears analogous to \cite{satia1973markovian}, it emphasizes the difference between the entire distributions rather than individual probabilities.  Generally, norm ambiguity sets can be represented as follows.
\begin{definition}[Norm ambiguity sets]
Let $\|\cdot\|$ denotes a general norm function, such as $L_1$-norm, $L_{\infty}$-norm and $L_p$-norm ($1<p<\infty$).  For $s\in\mathcal{S},a\in\mathcal{A}$, a $(s,a)$-rectangular norm ambiguity set can be defined as
    \begin{equation}
       \mathcal{P}_{sa}\triangleq\left.\left\{\mathbf{p}_{sa}\in\Delta^{|\mathcal{S}|} \right| \|\hat{\mathbf{p}}_{sa}-\mathbf{p}_{sa}\|\leq \theta_{sa},~\theta_{sa}\geq 0 \right\},
    \end{equation} 
    where $\hat{\mathbf{p}}_{sa}$ is the nominal distribution over the next state concerning state-action pair $(s, a)$, and $\mathbf{p}_{sa}$ belongs to $|\mathcal{S}|$-dimension probability simplex. The deviation with respect to the predetermined norm is constrained to be no larger than a parameter $\theta_{sa}$.
    For $s\in\mathcal{S}$, a $s$-rectangular norm ambiguity set can be defined similarly as
    \begin{equation}
        \mathcal{P}_s\triangleq\left\{\mathbf{p}_s=(\mathbf{p}_{s1},...,\mathbf{p}_{s|\mathcal{A}|})\in\Delta^{|\mathcal{S}|\times|\mathcal{A}|}\left|
        \begin{array}{l}\mathbf{p}_{sa} \in \Delta^{|\mathcal{S}|},\ \forall a\in \mathcal{A},\\
        \textstyle\sum\limits_{a\in\mathcal{A}}\|\hat{\mathbf{p}}_{sa}-\mathbf{p}_{sa}\|\leq \theta_s,\ \theta_{s}\geq 0\end{array}\right\}\right.,
    \end{equation}
    where $\mathbf{p}_{sa} \in \Delta^{|\mathcal{S}|},\ \forall a\in \mathcal{A}$ ensures that $\mathbf{p}_{sa}$ is a valid probability distribution over the next state, while  $\textstyle\sum_{a\in\mathcal{A}}\|\hat{\mathbf{p}}_{sa}-\mathbf{p}_{sa}\|\leq \theta_s$ establishes that the sum of the distances between the nominal probabilities $\hat{\mathbf{p}}_{sa}$ and plausible probabilities $\mathbf{p}_{sa}$ across all actions is bounded by $\theta_s$.
\end{definition}

Common norms, such as the $L_1$, $L_{\infty}$, and $L_p$-norm ($1<p<\infty$), are known to be convex. The convexity allows RMDPs with norm-based ambiguity sets to be readily formulated as convex programs, making them a preferred choice for tackling large-scale problems. Moreover, these norms often have intuitive geometric interpretations, which assist the decision-maker in establishing a desirable model. Among the various norms, the $L_1$ and $L_{\infty}$ norms have received significant attention in the literature.

The $L_1$-norm, also known as the \textit{variation distance}, is the most widely used and one of the earliest studied norms in the context of norm-constrained RMDPs \citep{iyengarRobustDynamicProgramming2005,petrik2014raam}. The construction of $L_1$-norm-based ambiguity sets offers two significant advantages. First, this construction enables the calculation of worst-case transition probabilities to be computed via linear programs, which brings computational efficiency. Second, the size of these ambiguity sets can be determined using Hoeffding-style bounds, which limit the probability of large deviations in the sums of random variables. This property is formalized in the following proposition.

\begin{proposition}[$L_1$-Norm Finite Sample Bound \citepalias{weissman2003inequalities}]
    \label{l1bound}
    Let $n_{sa}$ be the number of transitions from state $s$ by taking action $a$ and $\delta \in (0,1]$ be the confidence level. When the deviation bound parameter $\theta_{sa}$ is chosen as $\theta_{sa} = \sqrt{\frac2{n_{sa}}\log\frac{|\mathcal{S}||\mathcal{A}|2^{|\mathcal{S}|}}\delta}$, the underlying $\mathbf{p}_{sa}^{0}$ is contained in the ambiguity set with probability $1-\delta$.
\end{proposition}
 
Note that the bound in Proposition \ref{l1bound} gets tighter as $n_{sa}$ increases, reflecting higher confidence with more observations. The bound also depends on the size of the state and action spaces, as well as the desired confidence level $\delta$, all of which contribute to the complexity of achieving an accurate approximation. This bound allows for the data-driven construction of ambiguity sets by providing a principled way to balance between robustness and conservatism based on the amount of available data.

\citetalias{ho2018fast} consider a $w$-weighted $L_1$-norm constrained RMDPs with $(s,a)$- and $s$-rectangular assumption, where the norm is defined as $\|x\|_{1,w}=\textstyle\sum_{i=1}^{n}w_i|x_i|$.  In this paper, the authors leverage the properties of $L_1$-norm to develop fast Bellman updates, rather than directly using LPs. Specifically, the authors start from the Q-function $q\in\mathbb{R}^{|\mathcal{S}|\times|\mathcal{A}|}$ of RMDPs where:
\begin{equation}
	\label{qsa_function}
	q_{sa}(\eta)=\min_{\mathbf{p}_{sa}\in\Delta^{|\mathcal S|}}\{R(s,a)+\gamma\cdot\mathbf{p}_{sa}^Tv:\|\mathbf{p}_{sa}-\hat{\mathbf{p}}_{sa}\|\le\eta\}
\end{equation}
Where $R(s, a)$ is an immediate reward unrelated to the next state, $v:=(V(1),..., V(|\mathcal{S}|))$ denotes the state value vector and $\eta$ is controlled parameter and is smaller than predetermined robustness budget $\theta_{sa}$. As the robust Bellman optimality equation holds under rectangular assumptions \citep{iyengarRobustDynamicProgramming2005}, the authors derive optimal Q-function $q^*$ in return. Benefiting from the structure of $L_1$-norm, problem \eqref{qsa_function} can be readily reformulated as a linear program:
\begin{equation}
\label{l1norm_lp}
	\begin{array}{lll}
	q_{sa}(\eta)=&\min\limits_{\mathbf{p}_{sa}\in\mathbb{R}^{|\mathcal{S}|},~ l\in\mathbb{R}^{|\mathcal{S}|}} &  \mathbf{p}_{sa}^T(R(s,a)\mathbf{1}+\gamma v)\\
	&\text{s.t.} & \mathbf{p}_{sa}-\hat{\mathbf{p}}_{sa}\leq l\\
	&&\hat{\mathbf{p}}_{sa}-\mathbf{p}_{sa}\leq l\\
	&&\mathbf{p}_{sa}\geq\mathbf{0}\\
	&&\mathbf{1}^\mathsf{T}\mathbf{p}_{sa}=1, \quad w^\mathsf{T}l=\eta
	\end{array}
\end{equation}
Given the LP in \eqref{l1norm_lp} and $(s,a)$-rectangular assumption, the authors develop a homotopy method which starts from a trivial feasible solution $\eta=0$ (i.e., $\mathbf{p}_{sa}=\hat{\mathbf{p}}_{sa}$) and then track the optimal solution $\mathbf{p}_{sa}$ as $\eta$ gradually increases. Since $q_{sa}(\eta)$ and $\mathbf{p}_{sa}$ are both piecewise linear in $\eta$, the proposed homotopy approach can be efficiently implemented by tracing the basic feasible solutions. For the $s$-rectangular case, the authors show that an equivalent reformulation of the robust Bellman update exists, as follows:
\begin{equation}
\label{ho_reformulation}
	V^*(s)=\max_{\pi\in \Delta^{|\mathcal{S}|}}\min_{\eta\in\mathbb{R}^{|\mathcal{A}|}}\left\{\sum_{a\in \mathcal{A}}\pi(s,a)q_{sa}(\eta_a):\sum_{a\in \mathcal{A}}\eta_a\le \theta_s\right\}\iff \min_{u_s\in\mathbb{R}}\left\{u_s:\sum_a q_{sa}^{-1}(u_s)\leq\theta_s\right\}
\end{equation}
where $q_{sa}^{-1}(u)=\min_{\mathbf{p}_{sa}\in\Delta^S}\left\{\|\mathbf{p}_{sa}-\hat{\mathbf{p}}_{sa}\|_{1,w_a}:r_a+\gamma \mathbf{p}_{sa}^\mathsf{T}v\leq u_s\right\}$. This reformulation is non-trivial; therefore, we provide the intuition behind \eqref{ho_reformulation} instead of a concrete proof here. For a given state $s$, with a limited robustness budget $\theta_{s}$, $q_{sa}^{-1}(u)$ can be interpreted as the minimum robustness budget assigned to action $a$ such that the value-to-go function does not exceed $u_s$. Thus, minimizing $u_s$ determines the worst-case transition probabilities $\mathbf{p}_{sa}$ that lead to the lowest value-to-go function, which is equivalent to optimizing the robust Bellman equation. Note that this reformulation simplifies the original update into a one-dimensional problem, and the authors propose a bisection algorithm atop the homotopy method to efficiently solve it.

Although the homotopy method for the $(s, a)$-rectangular case in \cite{ho2018fast} has accelerated the search for the optimal solution, a linear program has to be solved for each state and each step of the value or policy iteration, which is costly in large state space. To overcome this disadvantage, \citetalias{ho2021partial} propose a partial policy iteration (PPI) framework for $(s, a)$- or $s$-rectangular ambiguity sets. In contrast to general robust policy iteration procedures \citep{iyengarRobustDynamicProgramming2005, nilimRobustControlMarkov2005}, the PPI simplifies the policy evaluation step where it only solves a regular ordinary MDP (which is constructed from the corresponding RMDP), thus any advanced MDP algorithm can be applied directly, making the iteration procedures more efficiently. Combining PPI with the homotopy method or bisection method, the resulting algorithms achieve a significant speedup. Additionally, the PPI proposed in this paper is the first policy iteration method that provably converges to the optimal solution for s-rectangular RMDPs, extending the previous findings that PPI merely holds for $(s, a)$-rectangular RMDPs \citep{kaufman2013robust}.

The $L_{\infty}$-Norm offers an alternative approach for constructing norm-based ambiguity sets, aiming to constrain the maximum deviation of each transition probability mass. Although $L_{\infty}$-norm is more intuitive and interpretable,
%closely resembling parametric approaches such as MDPIPs in section \ref{sec3}, 
and the corresponding ambiguity sets empirically outperform $L_1$-norm-based ambiguity sets in some circumstances \citepalias{pmlr-v130-behzadian21a}, it presents certain challenges compared to the $L_1$-norm. First, the concentration inequalities that facilitate the data-driven construction of high-confidence RMDPs do not hold. Second, most efficient algorithms, such as those in \cite{ho2018fast, ho2021partial}, rely on the sparsity properties of $L_1$-norm and therefore cannot be directly applied to $L_\infty$-norm ambiguity sets.

Inspired by \cite{ho2018fast, ho2021partial}, \citetalias{behzadian2021fast} also employ homotopy and bisection methods to design quasi-linear time complexity algorithms for $L_{\infty}$-constrained $(s, a)$ and $s$-rectangular RMDPs. The key procedure involves reformulating the Q-function as a parametric LP similar to \eqref{l1norm_lp}, but replacing the bounded constraints with
\begin{equation*}
	\mathbf{1}^T\mathbf{p}_sa=1,-\eta\le p_{sas'}-\hat{p}_{sas'}\le \eta,\ p_{sas}\ge 0, \ \forall s'\in\mathcal{S}.
\end{equation*}
Hence, leveraging the structural properties of the new LPs, such as the piecewise linear and non-increasing property of $q_{sa}(\eta)$ with respect to $\eta$, the authors modify the homotopy and bisection methods to accommodate the constraints brought by $L_\infty$-norm.

An interesting result for the choice between $L_1$-norm and $L_{\infty}$-norm is shown in \citetalias{russel2019optimizing}. From the perspective of value functions, they claim that the choice of set shape is dominantly driven by the structure of the value function, e.g., an $L_{\infty}$-constrained set is likely to work better than the $L_1$-constrained set when the value function is sparse. 

\subsection{\texorpdfstring{\(\phi\)}{phi}-Divergence Ambiguity}
While norm ambiguity sets offer advantages in computational efficiency and modeling convenience, they are limited by their ability to convey probabilistic information and guarantees. As a result, there is growing interest in ambiguity sets designed with probabilistic metrics. Among these, $\phi$-divergence stands out as a widely favored class due to its tractable and desirable statistical characteristics. Here, we formally define general $\phi$-divergence ambiguity sets.

\begin{definition}[$\phi$-divergence ambiguity] A $\phi$-divergence ambiguity set is defined via
    \begin{equation}
        \mathcal{P}\triangleq\left\{\mathbb{P}\in \mathfrak{M}(\Xi,\mathcal{F})\left|~ D^{\phi}(\mathbb{P},\hat{\mathbb{P}})\leq \theta\right\}\right.,
    \end{equation}
    where $D^{\phi}(\mathbb{P},\hat{\mathbb{P}}):=\int_{\Xi}\phi(\frac{\mathrm{d}\mathbb{P}}{\mathrm{d}\hat{\mathbb{P}}})\mathrm{d}\hat{\mathbb{P}}$ is the similarity measure function whose concrete representation depends on divergence function $\phi(\cdot):~\mathbb{R}_{+}\rightarrow \mathbb{R}_{+}\cup\{+\infty\}$. The divergence function $\phi(\cdot)$ is required to be convex, and satisfies $\phi(1)=0$, $0\cdot\phi(\frac{0}{0}):=0$ and $0\cdot \phi(\frac{a}{0}):=a\lim_{t\rightarrow\infty}\frac{\phi(t)}{t}$ for $a>0$ \citep{bayraksan2015data, pardo2018statistical}. 
\end{definition}

Despite numerous $\phi$-divergences, only some divergences are widely accepted due to their well-defined mathematical properties for robust optimization. Table \ref{phi-divergence functions} presents a list of typical divergence functions in terms of probabilities, where the second column $\phi(\cdot)$ is the corresponding $\phi$-divergence function, the third column $D^{\phi}(\mathbf{p}_{sa}, \mathbf{q}_{sa})$ is the representation of $\phi$-divergence in terms of probabilities, and last column is the conjugate function of $\phi(\cdot)$.

\begin{table}[htbp]
    \caption{Typical $\phi$-divergence functions and their conjugates $\phi^{*}(a)$ \citep{rahimian2019distributionally}}
    \label{phi-divergence functions}
\renewcommand\arraystretch{2}
    \centering
    \setlength{\tabcolsep}{2.7mm}{
    \begin{tabular}{llll}
    \toprule
    Divergence  & $\phi(t), t \geq 0$ & $D^{\phi}(\mathbf{p}_{sa}, \mathbf{q}_{sa})$ & $\phi^*(a)$ \\
    % \hline Kullback-Leibler  & $t \log t-t+1$ & $\sum_{s^{\prime}} p_{sas^{\prime}} \log(p_{sas^{\prime}} / q_{sas^{\prime}})$ & $e^a-1$ \\
    \midrule Kullback-Leibler  & $t \log t-t+1$ & $\sum_{s^{\prime}} p_{sas^{\prime}} \log(\frac{p_{sas^{\prime}}}{q_{sas^{\prime}}})$ & $e^a-1$ \\
     % Burg entropy  & $-\log t+t-1$ & $\sum q_\omega \log(q_{\omega} / p_{\omega})$ & $-\log (1-a), a<1$ \\
    Burg entropy  & $-\log t+t-1$ & $\sum_{s^{\prime}} q_{sas^{\prime}} \log(\frac{q_{sas^{\prime}}}{p_{sas^{\prime}}})$ & $-\log (1-a), a<1$ \\
    $J$-divergence  & $(t-1) \log t$ & $\sum_{s^{\prime}}\left(p_{sas^{\prime}}-q_{sas^{\prime}}\right) \log(\frac{p_{sas^{\prime}}}{q_{sas^{\prime}}})$ & No closed form \\
    $\chi^2$-distance  & $(t-1)^2/t$ & $\sum_{s^{\prime}} \frac{(p_{sas^{\prime}}-q_{sas^{\prime}})^2}{p_{sas^{\prime}}}$ & $2-2 \sqrt{1-a}, a<1$ \\
     Modified $\chi^2$-distance & $(t-1)^2$ & $\sum_{s^{\prime}} \frac{(p_{sas^{\prime}}-q_{sas^{\prime}})^2}{q_{sas^{\prime}}}$  & $\begin{cases} -1 & a<-2 \\
    a+\frac{a^2}{4} & a \geq-2\end{cases}$ \\
     Variation distance & $|t-1|$ & $\sum_{s^{\prime}}\left|p_{sas^{\prime}}-q_{sas^{\prime}}\right|$ & $\begin{cases} -1 & a \leq-1, \\
    a & -1 \leq a \leq 1\end{cases}$ \\
     Hellinger distance  & $(\sqrt{t}-1)^2$ & $\sum_{s^{\prime}}\left(\sqrt{p_{sas^{\prime}}}-\sqrt{q_{sas^{\prime}}}\right)^2$ & $\frac{a}{1-a},a<1$\\
    \bottomrule
    \end{tabular}}
\end{table}

Kullback-Leibler (KL) divergence, known for its interpretability and convexity, is among the most widely used $\phi$-divergences in both DRO and RMDPs. As the cornerstone of modern RMDPs, \cite{iyengarRobustDynamicProgramming2005} and \cite{nilimRobustControlMarkov2005} both improve upon traditional MDPIPs by exploring KL-divergence ambiguity sets to depict statistical uncertainty better. While sharing the same support between estimates and true randomness, KL-divergence limits the deviation in terms of the entire transition probability distribution, rather than individual probability mass. With $(s, a)$-rectangular assumption, they reformulate the inner maximization problem in dual form by employing the standard Lagrangian duality, as follows:
\begin{equation}
\label{valuefunction_reformulation}
\min_{\lambda>0}~\lambda\log\left(\sum\limits_{s'}\hat{p}_{sas^{\prime}}\exp\left(\frac{V_t(s'|s,a)}{\lambda}\right)\right)+\theta\lambda,
\end{equation}
where $\lambda$ is the Lagrangian multiplier and $\hat{p}_{sas^{\prime}}$ is the nominal probability of transition $(s,a,s^{\prime})$. As the problem is reduced to a one-dimensional convex program with respect to $\lambda$, \cite{nilimRobustControlMarkov2005} develop a bisection algorithm to obtain an $\epsilon$-optimal policy solved in polynomial time. Nonetheless, the reformulation problem still involves $\log$ and $\exp$ operators, posing certain computational difficulties.  To refine the theoretical foundations, \cite{iyengarRobustDynamicProgramming2005} utilizes the fact that $\log(1+x)\le x$ for all $x\in\mathbb{R}$, and proposes a conservative approximation ambiguity set in terms of probability measure,
\begin{equation}
\label{conservative_approximation}
    \mathcal{P}\triangleq\left\{P_{sa}\in \mathfrak{M}(\mathcal{S},\mathcal{F}^{\mathcal{S}})\ \middle|\ \sum\limits_{s'\in \mathcal{S}}\frac{(P_{sa}(s')-\hat{P}_{sa}(s'))^2}{\hat{P}_{sa}(s')}\leq \theta\right\} ,
\end{equation}
where $\mathcal{F}^{\mathcal{S}}$ is the $\sigma$-algebra associated with sample space $\mathcal{S}$ and the probability measure $P$ belongs to probability space $\mathfrak{M}(\mathcal{S},\mathcal{F}^{\mathcal{S}})$ that assigns probabilities to next state $s'$ occurring within a single period transition $(s,a,s')$ (with the time subscript omitted). Essentially, the summation in \eqref{conservative_approximation} represents a modified $\chi^2$-distance between $P$ and $\hat{P}$, the nominal probability measure. The author shows that this approximation allows us to exactly solve the inner optimization problem with the complexity $\mathcal{O}(|\mathcal{S}|\log|\mathcal{S}|)$, often more efficiently than original KL-divergence ambiguity set with complexity $\mathcal{O}(|\mathcal{S}|^{1.5}\log|V_{max}/\epsilon|)$ where the solution is $\epsilon$-optimal.

Motivated by the findings in \cite{iyengarRobustDynamicProgramming2005} and \cite{nilimRobustControlMarkov2005}, \cite{liu2022distributionally} design a distributionally robust Q-learning algorithm with the KL-divergence ambiguity sets. Like \eqref{valuefunction_reformulation}, the authors leverage the existence of the robust Bellman equation (arising from $(s, a)$-rectangular assumption) and the strong duality lemma from classical DRO results under KL-perturbation \citep{hu2013kullback} to derive a distributionally robust Q-function, which transforms the primal infimum (inner) infinite-dimensional problem into a supremum finite-dimensional dual representation\footnote{For simplicity, we demonstrate the distributionally robust Q-function using a simplified formulation without reward uncertainty that is slightly modified from the original version of the paper.}:
\begin{equation}
    \begin{aligned}
        Q(s,a):=R(s,a)+\gamma \sup_{\lambda \geq 0}\left\{- \lambda\log\left(\mathbb{E}_{\hat{\mathbf{p}}_{s,a}}\left[\exp\left(-\frac{\max_{b\in\mathcal{A}}Q(s',b)}{\lambda}\right)\right]\right)-\lambda\theta\right\}.
    \end{aligned}
\end{equation}
where $R(s,a)$ is the immediate reward unrelated to the next state and $\gamma$ is the discount factor. The supremum is taken over the non-negative dual variable $\lambda$ derived from the primal problem transformation. The key term $\lambda\log(\mathbb{E}[\cdot])$ results from the dual formulation of the distributionally robust problem, where the expectation taken over the nominal transition probabilities  $\hat{\mathbf{p}}_{s, a}$, ensuring that the reformulation becomes finite-dimensional and tractable. In essence, the worst-case Q-function is equivalent to a plug-in estimator over nominal transition probabilities when the deviation measured by the KL-divergence is not larger than $\theta$.

While the (dual) distributionally robust Bellman equation is well-established, obtaining an unbiased estimator using a simulator that samples from $\hat{\mathbf{p}}_{sa}$ remains challenging due to its nonlinear structure. This nonlinearity stems from $\log$ function and the nested expectation term. In such nonlinear settings, simply substituting sample average approximations for true expectations does not yield unbiased estimates, as the expectation of a nonlinear function of random variables generally differs from the function of the expectations of those variables. To address this, the authors introduce the multi-level Monte-Carlo scheme \citepalias{blanchet2015unbiased, blanchet2019unbiased}, which reduces the final estimation bias through introducing $N$ estimators with varying degrees (i.e., increasing precise degrees). With the new unbiased robust Q-value estimates derived from the simulator samples, the authors demonstrate that the proposed algorithm converges asymptotically to the optimal distributionally robust problem. Notably, to the best of our knowledge, this proposed algorithm is the first model-free algorithm ever developed on RMDPs.

Another commonly used approach to model the ambiguity sets in RMDPs is $\chi^2$-divergence. Compared to KL-divergence, $\chi^2$-divergence possesses several advantageous properties, including symmetry, reduced sensitivity to minor differences, and computational convenience (recall the conservative approximation proposed by \cite{iyengarRobustDynamicProgramming2005} above is essentially $\chi^2$-divergence).  \cite{hanasusanto2013robust} consider a data-driven stochastic control problem with continuous state and action spaces, where the transition probability matrix is estimated via Nadaraya-Watson (NW) kernel regression, as the nominal distribution. Given finite horizon $T$ and $N$ sample trajectories $\{\xi^{i}_{t}\}^{T}_{t=1},~i\in[N]$, the authors use empirical estimates to approximate the conditional expected value of future state $V_{t+1}$ concerning current observation $\xi_t$, as follows:
\begin{equation}
\label{nw_approximation}
    \mathbb{E}[V_{t+1}(s_{t+1},{\xi}_{t+1})~|~{\xi}_t]\approx\sum_{i=1}^{N} p_{ti}({\xi}_t)V_{t+1}(s_{t+1}^i,{\xi}_{t+1}^i),
\end{equation}
where $p_{ti}({\xi}_t)=\frac{\mathrm{K_{NW}(\xi_t-\xi^{i}_t)}}{\sum_{j=1}^{N}\mathrm{K_{NW}}(\xi_t-\xi^{j}_t)},~i\in[N],t\in[T]$ and $\mathrm{K_{NW}}(\cdot)$ is the NW kernel density function. In the right-hand side of \eqref{nw_approximation}, the weights give more importance to samples that are ``closer" to the current observation $\xi_t$ in the future space. This method provides a non-parametric way to estimate conditional expectations, which is particularly useful in continuous state spaces where traditional tabular methods might not be applicable. The NW kernel regression allows for smooth interpolation between observed data points, providing estimates even for states not directly observed in the training data.

However, if the training data is sparse, the NW estimates typically exhibit high variance, leading to the poor out-of-sample performance of the data-driven DP. To mitigate the impact of this issue, the authors construct $\chi^2$-distance ambiguity sets. Implicitly assuming $(s, a)$-rectangular assumption, they reformulate the inner infinite-dimensional maximum problem as a tractable minimum problem by standard dual theory. By exploiting the structure of $\chi^2$-distance ambiguity sets, the reformulation can be expressed as a second-order cone program (SOCP) under specific conditions. For instance, if the value function is piecewise linear or convex quadratic, and other mild assumptions hold—such as the immediate reward or cost function being convex quadratic in the state and action—this SOCP can be efficiently solved using interior-point algorithms. Instead of finding the optimal policies, they design a robust data-driven dynamic program algorithm to approximate the value-to-go functions via interpolation.

Also adhering to $\chi^2$-divergence, \citetalias{klabjan2013robust} directly utilize historical data to construct ambiguity sets from a goodness-of-fit test perspective for the single-item multi-period periodic review stochastic lot-sizing problem. Different from the approach in \cite{hanasusanto2013robust}, the authors construct the nominal distribution discretely. Let $B$ represent the set of bins that partition the range of possible values of the random variable $\xi_t\sim \mu_t$. For each $i\in B$, let $N_{t,i}(\xi_t)$ be the number of observations falling within the $i$-th bin, %Let the number of observations falling within the \textit{i}-th bin be a function of a realized random variable $\xi_t\sim \mu_t$ as $N_{t,i}(\xi_t)$, 
and the total number $n_t(\xi_t)=\sum_{i\in B}N_{t,i}(\xi_t)$. The ambiguity set consists of the distributions $\mu_t$ which satisfy
\begin{equation}
\label{hypothesis_test}
    \sum_{i\in B}\frac{(N_{t,i}(\xi_t)-n_t(\xi_t)\cdot\hat{\mu}_t(i))^2}{n_t(\xi_t)\hat{\mu}_t(i)}\leq \theta,\quad t=1,...,T,
\end{equation}
where $\hat{\mu}_t(i)$ denotes the estimated probability that an observation falls within the $i$-th bin. The parameter $\theta$ controls how close the observed sample data is to the estimated expected number of observations according to the fitted distribution $\mu_t$. From a hypothesis testing perspective, the authors set $\theta=\chi^2_{|B|-1,1-\alpha}$ where $|B|$ represents the total number of bins and $\alpha$ denotes the significance level. Consequently, all distributions satisfying the constraints in \eqref{hypothesis_test} are those for which the corresponding null hypothesis is not rejected at the significance level $\alpha$. By implementing a standard dualization process, the original robust problem can be converted into a tractable second-order cone program. In this paper, the authors show that a state-dependent $(s, S)$ policy is optimal for this robust lot-sizing problem and the ordering levels can be computed by this second-order cone program.

Notably, recent literature on RMDPs with $\phi$-divergence ambiguity sets has increasingly focused on sample complexity and the design of efficient algorithms. Unlike the efficient algorithms for RMDPs with norm ambiguity sets, which typically necessitate the use of LPs, the efficient algorithms involving $\phi$-divergence display distinct methodologies. In these aspects, there have been several noteworthy and impressive advancements. Below, we briefly review some of the advanced literature in this area.

For sample complexity analysis, \citetalias{wang2023finite} extend the distributionally robust Q-learning framework proposed by \cite{liu2022distributionally} through refining the design and analysis of the key component, multi-level Monte Carlo estimator. Specifically, while the expected number of samples requested in \cite{liu2022distributionally} is infinite and the algorithm only can achieve asymptotic guarantees, the refined algorithm in \citetalias{wang2023finite} requires a merely constant order number of samples. This improvement is achieved by delicately devising the sample numbers for each state-action pair. This is also the first sample complexity result for the model-free robust RL problem. \citetalias{shi2024curious} investigate RMDPs with the uncertainty set measured via total variation (TV) distance or $\chi^2$-divergence. The authors suppose that one has access to a generative model or simulator to draw samples with a nominal transition probability matrix. Under $(s, a)$-rectangular assumption, they propose a model-based distributionally robust value iteration algorithm. In each iteration, the robust Q-function is computed in terms of each state-action pair, and the robust value function is set greedily according to the current robust Q-function. Interestingly, the authors find that the choice of uncertainty set significantly impacts the sample size requirements. For instance, to achieve the same $\epsilon$-accuracy, while the RMDPs with $\chi^2$-divergence are harder than standard MDPs as expected, the RMDPs with TV distance are easier than standard MDPs. These findings emphasize the importance of the construction of uncertainty sets and the efficiency of the model-based approaches.
%Meanwhile, they strengthen the lower and upper sample complexity bound for both metrics.} 

For efficient algorithms design, \cite{grand2021scalable} introduce the first first-order method (FOM) for solving RMDPs and propose a scalable algorithmic framework of FOM updates with occasional approximate value iteration updates (FOM-VI) to limit the computational cost of value iteration. This framework represents the first tractable algorithm for $s$-rectangular RMDPs with KL-divergence ambiguity sets. It is noteworthy that previously, only the $(s, a)$-rectangular case had been effectively addressed by robust value iteration in \cite{nilimRobustControlMarkov2005}.
%In summary, the FOM-VI framework interleaves FOM updates with occasional approximate VI updates, thereby effectively limiting the computational cost of VI updates. 
A key observation underpinning this algorithm is that an $s$-rectangular RMDP can be decomposed into $|\mathcal{S}|$ bilinear saddle-points problems (BSPPs). Within each epoch, a primal-dual algorithm (PDA) is employed to update the action policy $\pi_s$ (primal updates) and the conditional transition matrix $\mathbf{p}_s$ (dual updates) with fixed value-to-go function $V_s$ for each state $s$ (i.e., a BSPP). At the end of this epoch, $V_s$ is updated by fixed policy and conditional transition matrix, which are constructed by a weighted average of the policies or conditional transition matrices with designed weights.  A novel feature of this algorithm lies in the PDA updates, where the associated proximal mappings are defined as follows:
\begin{equation}
\label{prox_mappings}
	\begin{aligned}
		\mathrm{prox}_{\pi}(\boldsymbol{g}_{\pi},\pi_s')&=\arg\min_{\pi_s\in\Delta^{|\mathcal{A}|}}\langle\boldsymbol{g}_{\pi},\pi\rangle+D_{\pi}(\pi_s,\pi_s^{\prime}),\\
		\mathrm{prox}_{\mathbf{p}}(\boldsymbol{g}_{\mathbf{p}},\mathbf{p}_s^{\prime})&=\arg\max_{\mathbf{p}_s\in\mathcal{P}_s}\langle\boldsymbol{g}_{\mathbf{p}},\mathbf{p}_s\rangle-D_{\mathbf{p}}(\mathbf{p}_s,\mathbf{p}_s^{\prime}),
	\end{aligned}
\end{equation}
where $\boldsymbol{g}_{\pi}$ and $\boldsymbol{g}_{\mathbf{p}}$ are the gradients with respect to current $\pi_s$ and $\mathbf{p}_s$, respectively, and $D_{\pi}(\cdot,\cdot)$ and $D_{\mathbf{p}}(\cdot,\cdot)$ are Bregman divergence functions. The intuition behind \eqref{prox_mappings} is that proximal mappings move along with the direction of improvement as indicated by the gradients while being penalized by the Bregman divergences. These penalties ensure that the updates remain within a region where the first-order approximations are accurate enough. Given the KL-divergence constraints in $\mathcal{P}_s$, the authors show that there exists a Bregman divergence such that solving the original RMDPs is equivalent to operating the FOM-VI framework, where the PDA updates are governed by the proximal mappings in \eqref{prox_mappings}.

\cite{ho2022robust} further generalize existing literature by proposing a fast algorithmic framework to solve general $\phi$-divergence RMDPs for both $(s, a)$-rectangular and $s$-rectangular ambiguity sets. Unlike \cite{grand2021scalable}, which replaces traditional VI procedures with FOMs updates, the authors adhere to previous research methodologies in RMDPs, focusing on accelerating the robust Bellman update in robust VI. The key challenge in solving $s$-rectangular RMDPs with $\phi$-divergence ambiguity sets lies in the fact that the corresponding ambiguity sets are non-polyhedral. The non-polyhedral nature makes each iteration (i.e., the robust Bellman update) computationally costly, a problem not present in the $(s, a)$-rectangular counterparts. This explains why the $(s, a)$-rectangular case can be efficiently solved by robust VI, whereas the $s$-rectangular case cannot. 
Given this barrier, a crucial component in this paper that circumvents the aforementioned challenge is the authors' discovery that these min-max problems can be reduced to a small number of highly structured projection problems onto a probability simplex.

Specifically, by applying the minimax theorem, the optimal value of the original max-min problem equals to min-max problem where the inner problem is optimized over a deterministic worst action, rather than a (randomized) policy. The min-max problem can be efficiently solved via the bisection method; which finds the lowest possible constant $\beta$ such that the objective value is not larger than $\beta$ for any feasible solution. Once $\beta$ is given by the bisection method, the following generalized $D_{\phi}$-projection problem of nominal transition probabilities $\hat{\mathbf{p}}_{sa}$ will be checked:
\begin{equation}
\label{projection_problem}
	\begin{array}{llll}
	\mathrm{P}(\hat{\mathbf{p}}_{sa};\boldsymbol{b},\beta)= &\mathrm{min}& &D_{\phi}(\mathbf{p}_{sa},\hat{\mathbf{p}}_{sa})\\
	&\mathrm{s.t.}& &\boldsymbol{b}^\top\mathbf{p}_{sa}\leq\beta\\
	&& &\mathbf{p}_{sa}\in\Delta^{|\mathcal{S}|}
	\end{array}
\end{equation}
where $D_{\phi}$ is the predetermined $\phi$-divergence function and $\boldsymbol{b}\in\mathbb{R}^{|\mathcal{S}|}_+$ is a known parameter vector. Set $\boldsymbol{b}=\boldsymbol{r}_{sa}+\lambda \boldsymbol{v}$ and compute the value $\sum_{a\in\mathcal{A}}\mathrm{P}(\hat{\mathbf{p}}_{sa};\boldsymbol{r}_{sa}+\lambda \boldsymbol{v},\beta)$ and the corresponding optimal solution $\mathbf{p}^*_{sa}$. Then, the authors distinguish between the following two cases:
\begin{enumerate}
	\item If $\sum_{a\in\mathcal{A}}\mathrm{P}(\hat{\mathbf{p}}_{sa};\boldsymbol{r}_{sa}+\lambda \boldsymbol{v},\beta)\le \theta$, then $\mathbf{p}_{s}=(\mathbf{p}^*_{sa})_{a\in\mathcal{A}}$  is a feasible solution to \eqref{projection_problem}. Consequently, $\beta$ upper bounds the optimal objective value of the original RMDP.
	\item If $\sum_{a\in\mathcal{A}}\mathrm{P}(\hat{\mathbf{p}}_{sa};\boldsymbol{r}_{sa}+\lambda \boldsymbol{v},\beta)> \theta$, then there is no feasible $\mathbf{p}_{s}\in \Delta^{|\mathcal{S}|\times|\mathcal{A}|}$ such that the  objective value attains $\beta$ or less. Thus, $\beta$ becomes a lower bound of the optimal objective value.
\end{enumerate}

As $\mathrm{P}(\hat{\mathbf{p}}_{sa};\boldsymbol{r}_{sa}+\lambda \boldsymbol{v},\beta)$ represents the lowest deviation between $\hat{\mathbf{p}}_{sa}$ and $\mathbf{p}_{sa}$ that satisfies the value-go-function is bounded by $\beta$ for each state-action pair $(s,a)$, the $\sum_{a\in\mathcal{A}}\mathrm{P}(\hat{\mathbf{p}}_{sa};\boldsymbol{r}_{sa}+\lambda \boldsymbol{v},\beta)$ is the total deviation between $\hat{\mathbf{p}}_{s}$ and $\mathbf{p}_{s}$. If case 1 holds, the optimal solution $\mathbf{p}_{s}$ is a feasible solution in which each $\mathbf{p}_{sa}$ is the worst-case. By duality theory, $\beta$ becomes the upper bound of the optimal value. In turn, if the deviation $\sum_{a\in\mathcal{A}}\mathrm{P}(\hat{\mathbf{p}}_{sa};\boldsymbol{r}_{sa}+\lambda \boldsymbol{v},\beta)$ exceeds the parameter $\theta$, the solution $\mathbf{p}_{s}$ is out of the constructed ambiguity set, and corresponding $\beta$ provides a lower bound. Given this analysis and relationship, consequently, if the projection problem like \eqref{projection_problem} can be solved efficiently, the robust Bellman update of RMDPs can also be computed efficiently by the bisection method.

\subsection{Wasserstein Distance Ambiguity}
While $\phi$-divergence has demonstrated strong performance in certain applications, it also has notable limitations that hinder its modeling flexibility. For instance, $\phi$-divergence operates only on distributions with the same support and does not satisfy the triangle inequality. In contrast, Wasserstein distance has emerged as a good substitute for $\phi$-divergence. Wasserstein distance is a particular case of optimal transport discrepancies, which computes the minimal cost of transporting the masses between two distributions. Based on the definition of Wasserstein distance, the Wasserstein distance ambiguity is defined as follows.

\begin{definition}[Wasserstein Distance Ambiguity] The Wasserstein distance $\mathcal{W}(\mathbb{P},\hat{\mathbb{P}})$ of $\mathbb{P},\hat{\mathbb{P}}\in \mathfrak{M}(\Xi,\mathcal{F})$ is defined via
    \begin{equation}
    \label{wasserstein_distance}
        \mathcal{W}(\mathbb{P},\hat{\mathbb{P}})=\min\limits_{\Gamma\in\mathfrak{M}\times\mathfrak{M}}\left\{\int_{\Xi\times \Xi}c(\xi, \zeta)\Gamma(\mathrm{d}\xi,\mathrm{d}\zeta):~ \Pi^{1}_{\#}\Gamma=\mathbb{P},~\Pi^{2}_{\#}\Gamma=\hat{\mathbb{P}}\right\},
    \end{equation}
    where $c(\xi, \zeta)$ is the symmetric cost function of moving the mass from $\xi$ to $\zeta$ (i.e., $c(\xi,\zeta)=c(\zeta,\xi)$), $\Pi^{i}_{\#}\Gamma$ denote the $i$-th marginal distribution of $\Gamma$, and $\Gamma(\cdot,\cdot)$ is the joint distribution of $\mathbb{P}$ and $\hat{\mathbb{P}}$, also called the transport plan. Consequently, given the radius $\theta$ that controls the level of robustness, a Wasserstein ambiguity set can be constructed as
    \begin{equation}
        \mathcal{P}\triangleq\left\{\mathbb{P}\in \mathfrak{M}(\Xi,\mathcal{F})\ \middle|\ \mathcal{W}(\mathbb{P},\hat{\mathbb{P}})\leq \theta\right\}.
    \end{equation} 
\end{definition}  

While the Wasserstein distance can be defined well in terms of a joint distribution $\Gamma$ in \eqref{wasserstein_distance}, its dual form, also referred to as the Kantorovich metric, sometimes contributes to the modeling and analysis. Formally, we claim that the following result is equivalent to the primal definition of Wasserstein distance in \eqref{wasserstein_distance}.
    
\begin{proposition} Wasserstein distance has a dual representation due to Kantorovich duality \citep{villani2009optimal}, and we can represent Wasserstein distance in \eqref{wasserstein_distance} equivalently as
\begin{equation}
\label{Kantorovich}
    W(\mathbb{P},\hat{\mathbb{P}}) = \sup_{\substack{\psi\in L^{1}(\mathbb{P})\\\varphi\in L^{1}(\hat{\mathbb{P}})}}\left\{\int_{\Xi}\psi(\xi)\mathbb{P}(\mathrm{d}\xi)+\int_{\Xi}\varphi(\zeta)\hat{\mathbb{P}}(\mathrm{d}\zeta):~\psi(\xi)+\varphi(\zeta)\leq c(\xi,\zeta)~ \forall \xi,\zeta\in \Xi\right\},
\end{equation}
where $L^{1}(\mathbb{P})$ denotes the $L^1$ space of functions that are $\mathbb{P}$-measurable. Especially, when the $c(\cdot,\cdot)$ is 1-order, namely, 1-Wasserstein distance is considered, the dual form can be further simplified as
\begin{equation}
	\label{1_wasserstein_dual}
	W(\mathbb{P},\hat{\mathbb{P}}) = \sup_{\psi\in L^{1}(\mathbb{P})}\left\{\int_{\Xi}\psi(\xi)\mathbb{P}(\mathrm{d}\xi)-\int_{\Xi}\psi(\xi)\hat{\mathbb{P}}(\mathrm{d}\zeta)\right\},
\end{equation}
\end{proposition}

The dual formulation transforms the original minimization problem over transport plans $\Gamma$ into a maximization problem over potential functions $\psi,\varphi$,  subject to a point-wise constraint. The dual form simplifies the original complex problem with coupling constraints as an optimization over 1-Lipschitz functions. Additionally, the dual formulation provides a framework for establishing theoretical results, such as convergence rates.

Prominent results from recent Wasserstein DRO literature \citep{mohajerin2018data, blanchet2019quantifying, gao2023distributionally} offer tractable reformulation procedures and finite-sample guarantees, demonstrating that how to construct Wasserstein ambiguity sets in a data-driven manner and inspiring further research on RMDPs with Wasserstein ambiguity sets.

The Wasserstein DRMDPs were first investigated by \cite{yang2017convex}, where the elements of the ambiguity sets are the joint distribution of transition probabilities and immediate rewards. Implicitly assuming $s$-rectangularity, the author derives the corresponding distributionally robust Bellman equation by applying DP principles
\begin{equation}
V_t(s)=\sup_{\pi\in\Delta^{|\mathcal{A}_s|}}\inf_{\mu\in\mathcal{P}_s}\int_{\Xi_s}\sum_{a\in\mathcal{A}_s}\pi(s,a)\left(R_s(s,a)+\textstyle\sum\limits_{s^{\prime}\in\mathbb{S}}P_{s}(s^{\prime})V_{t+1}(s^{\prime})\right)\mathrm{d}\mu(P_s,R_s),
\end{equation}
and shows that this DRMDP admits an optimal Markov policy, while its construction is computationally challenging. To address the computational issue, the author rewrites the Wasserstein distance as the dual form and leverages Kantorovich duality \citep{villani2021topics} to obtain a dual reformulation. Using this reformulation, the finite horizon DRMDPs can be solved by finite-dimensional convex programming. Motivated by the findings in \cite{gao2023distributionally}, the author offers a closed form of the worst-case probability distribution when the nominal distribution has finite support. However, no efficient algorithm or supplementary statistical properties, such as finite-sample guarantees, associated with Wasserstein distance are provided in this paper.

Building on the foundation work of \cite{yang2017convex}, \cite{yang2020wasserstein} considers a distributionally robust optimal control problem and extends the finite state space setting of \cite{yang2017convex} into a continuous state space setting. Although optimal control problems slightly differ from MDPs, as previously mentioned, they can be connected under some mild settings. Consequently, employing similar procedures and methodologies, the author also shows the existence of an optimal policy and provides a tractable dual reformulation with no duality gap. Given the $s$-rectangular assumption, the author employs value iteration and policy iteration algorithms within this dual problem framework, both of which involve solving semi-infinite programs in each iteration. The author also establishes the least number of iterations required to attain an $\epsilon$-optimal policy. Meanwhile, the author tailors specific properties from Wasserstein DRO literature to characterize proposed Wasserstein RMDPs. Under certain mild assumptions, the optimal policy for the worst-case scenario is deterministic and stationary, and its structure can be expressed explicitly, as shown in \cite{gao2023distributionally}, where Kantorovich duality and DP play a critical role. Furthermore, a probabilistic out-of-sample performance guarantee with mild assumptions is established, similar to the results in \cite{mohajerin2018data}.

Apart from DRMDPs with Wasserstein distance, \cite{ramani2022robust} consider $(s, a)$-rectangular RMDPs with general distance metrics under both finite horizon and infinite horizon settings, thereby including Wasserstein distance, where the nominal transition probabilities are derived by sample average approximation (SAA). As a theoretical paper, \cite{ramani2022robust} aims to establish the following results: ($i$) robust value convergence, ($ii$) probabilistic performance guarantees on out-of-sample values, and ($iii$) a probabilistic convergence rate in the infinite horizon of rectangular RMDPs. For the first claim, the authors prove that the optimal values of the RMDPs converge almost surely to the true optimal value as the sample size $N\to\infty$. The key assumption required for the proof is the two probability mass functions are close if the metric function used in ambiguity sets deems them to be close. The intuition behind this assumption is that the simultaneous convergence of the empirical estimates to the true transition probabilities, along with the shrinking radius of the ambiguity balls, ensures convergence. Regarding the second claim, the rectangularity allows the probability of arbitrary guaranteed performance to be decomposed into the Cartesian product of the probabilities for each $(s, a)$ pair, facilitating further contractions. In the infinite horizon, by leveraging the \textit{simulation lemma} in the RL literature \citep{rajeswaran2020game}, the authors bound the difference between the values of two policies through the corresponding transition matrices and thereby derive the third claim.

The successful application of Wasserstein distance in deep learning and classical RMDPs has inspired a growing body of research that integrates it with reinforcement learning algorithms. However, a considerable portion of literature adopts a \textit{two-step} framework that treats Wasserstein distance as an independent plug-in component firstly, such as for updating simulator \citep{abdullah2019wasserstein} or restricting feasible action set \citep{kandel2020safe}, and then solves a modified non-robust MDP. Not surprisingly, this two-step framework lacks theoretical guarantees and statistical properties that should be derived from the Wasserstein distance.

With stage-wise independent assumption, implying the $(s, a)$-rectangularity holds, \cite{neufeld2024robust} design a robust Q-learning algorithm with Wasserstein ambiguity sets. Different from previous works that utilize Wasserstein distance separately, the authors leverage the dual results of the robust Bellman equation due to $(s, a)$-rectangularity, and reformulate the Q-value update in a regularized form. Specifically, they introduce a so-called $\lambda c$-transform of function $f$, which is essentially a simplified and straightforward outcome resulting from strong duality as
\begin{equation*}
	(f)^{\lambda c}(x):=\sup_{y\in\mathcal{X}}\{f(y)-\lambda\cdot c(x, y)\}
\end{equation*}
where $\mathcal{X}$ is the support of $x, y$, $c: \mathcal{X}\times \mathcal{X}\to \mathbb{R}$ is a predetermined cost function, and $f: \mathcal{X}\to \mathbb{R}$. In this paper, if we let $f$ be the value-to-go function and $c$ be the Wasserstein distance, the $\lambda c$-transform is the dualization procedure, omitting the Lagrangian multiplier term. This proposed algorithm is primarily rooted in the classical Q-learning algorithm. Nonetheless, before updating Q-values in each iteration, the optimal Lagrangian multiplier $\lambda^*$ will be computed such that the $\lambda c$-transformed Bellman equation with the Lagrangian multiplier term attains the supremum taken over the nominal probabilities:
\begin{equation}
    \label{optimal_lambda}
	\lambda^*=\arg \sup_{\lambda\geq0}\mathbb{E}_{\hat{\mathbf{p}}_{sa}}\left[-(-f_{(s,a)})^{\lambda c}(s')-\varepsilon\lambda\right], 
\end{equation}
where $f_{(s,a)}$ denotes the value-to-go function with respect to $(s,a)$. Subsequently, $\lambda^*$ will be used for the Q-values update where the value of dual form safely replaces the original (worst-case) Q-value, due to strong duality. The authors rigorously prove that the output of the proposed algorithm converges to the optimal robust Q-value function. Furthermore, the difference between the values of the robust and non-robust MDPs can be bounded if the Wasserstein ambiguity sets contain the true probability.

More generally, most RL settings lack a clear characterization of the environment, specifically in terms of access to a nominal distribution $\hat{\mathbf{p}}$. Consequently, model-based approaches, such as the one introduced by \cite{neufeld2024robust}, become inapplicable (see the optimization procedure in (42) for reference). To address this issue, \cite{hou2020robust} introduce a model-free algorithm called Wasserstein Robust Advantage Actor-Critic (WRAAC), which embeds Wasserstein constraints into standard actor-critic-based algorithms. The authors apply the Lagrangian method and the strong duality property brought by Wasserstein distance \citep{blanchet2019quantifying} to convert the min-max robust Bellman equation into a finite-dimensional optimization problem under an implicit $(s, a)$-rectangular assumption. As they demonstrate the existence of a deterministic Markov optimal policy and value-to-go function for the reformulated expected Bellman-form operator, they design the WRAAC algorithm where two neural networks are constructed to act as the critic and actor: the critic serves as the estimator of the value function, the actor aims to attain the robust optimal policy. WRAAC is implemented as a double-loop algorithm: the inner loop determines the extent of perturbations, and the outer loop optimizes the policy as normal procedures. In the inner loop, the algorithm updates the state $z$ that maximizes the penalized value-to-go function (i.e., by the Wasserstein distance with Lagrangian multiplier) and the Lagrangian multiplier $\lambda^*$ by sampling. In essence, the inner loop approximates the procedure \eqref{optimal_lambda} by Monte Carlo. With the currently optimal $\lambda^*$, the outer loop optimizes the critic and actor by the temporal errors. After the actor is updated, new samples are collected from the environment using the updated policy, and the process repeats. Sharing a similar idea to \cite{neufeld2024robust}, the multiplier $\lambda^*$ condenses the impact of the worst-case scenario, which is why we argue that the inner loop certifies the extent of perturbations.

\section{Beyond Rectangularity: Coupled Uncertainty Modeling}\label{sec6}

While rectangular assumptions bring tractability to RMDPs, such rectangular uncertainty sets are overly conservative in modeling uncertainty and are not always appropriate in practice \citep{iyengarRobustDynamicProgramming2005}. The following example demonstrates the necessity of relaxing the rectangular assumptions.
\begin{example}[Health Evolution, \cite{goyalRobustMarkovDecision2023}]     \label{nonrectangularexample}
    A Markov model may be used to describe the health evolution of a patient. The state $s\in \mathcal{S}$ represents the health condition, and the action $a\in A$ represents the treatment. The transition probabilities $\mathbf{p}_{sa}$ represent the dynamics of the health evolution across different health conditions given a treatment. In this context, factors such as genetics and disease traits can influence how patients transition between states, resulting in uncertainty correlation in the transition probabilities across different states. For example, it is reasonable to believe that two patients in the same moderate state may have different transition probabilities based on their previous conditions: one transitioning from mild to moderate (indicating a worsening trend) and the other from severe to moderate (indicating a recovery trend).
\end{example}

\citetalias{mannor2012lightning} first propose the coupled uncertainty sets LDST motivated by the famous proverb ``\textbf{L}ightning \textbf{D}oes not \textbf{S}trike \textbf{T}wice". As the proverb goes, the motivation behind LDST is that if the parameters of each state deviate with a small probability, and all states are independent, then the total number of states with deviated parameters will be small. Concretely, the uncertainty sets  consist of support constraints and cardinality constraints which limit the number of times parameters can deviate from the nominal values by a given parameter $K$, i.e., at most $K$ states can have transition probabilities different from their nominal estimates:
\begin{equation}
    \label{lightning}
    \mathcal{P}=\left\{\mathbf{p}\in (\Delta^{|\mathcal{S}|})^{|\mathcal{S}|\times |\mathcal{A}|}\left|\mathbf{p}_s\in \mathcal{P}_s,~\forall s\in \mathcal{S};~ \textstyle\sum\limits_{s\in \mathcal{S}}\boldsymbol{1}(\mathbf{p}_s\neq \hat{\mathbf{p}}_s)\leq K\right\}\right..
\end{equation}

The authors consider two models: a non-adaptive model, which can be viewed as a single-stage game over $T$ time steps yet with a joint cardinality constraint; and an adaptive model, where both the decision maker and the nature are aware of the number of ``deviated visits" over time and modify their strategies dynamically. While similar probabilistic guarantees are provided for both models, the authors demonstrate that solving the non-adaptive model is generally NP-hard, as shown through a reduction from the \textit{Vertex Cover Problem}. However, the non-adaptive model becomes tractable when limited to reward uncertainty.  In contrast, with the deviation information, adaptive models can be tractable by augmenting the state space as
\begin{equation}
    \overline{\mathcal{S}}_K=\mathcal{S}\times[0:K],~ \overline{\mathcal{S}}_N=\mathcal{S}\times[0:K]\times\mathcal{A},~\overline{\mathcal{S}}=\overline{\mathcal{S}}_K\bigcup\overline{\mathcal{S}}_N,
\end{equation}
where $\overline{\mathcal{S}}_K$ denotes the states of decision maker, and $\overline{\mathcal{S}}_N$ denotes the states of the nature. For each $s \in \mathcal{S}$, $k\in [0:K]$, the decision-maker's state $(s,k)$ denote the original state $s$ and deviation information $k$. Suppose the decision-maker takes action $a\in\mathcal{A}$, the next state will be $(s,k, a)$ for nature. With the augmented state space $\overline{\mathcal{S}}$ and the Nash equilibrium construction, the adaptive models are shown to be solved in polynomial time via backward induction or value iteration. Besides its tractability under non-rectangular settings, another key benefit of LDST is that it requires no distribution information.

Adhering to the same inspiration, \cite{mannor2016robust} further generalizes LDST ambiguity sets to the so-called $k$-rectangular ambiguity sets, treating the LDST ambiguity set as a special case. To ensure readability and completeness, we introduce the definition of conditional projection which is used for $k$-rectangular ambiguity sets. 
\begin{definition}[Conditional Projection]
    Let $\mathcal{S}'\subset \mathcal{S}$ be a nonempty subset of the states. For a subset $\mathcal{S}'\subset \mathcal{S}$, we denote $\mathbf{p}_{\mathcal{S}'}=\bigotimes_{s\in\mathcal{S}'}\mathbf{p}_s$. The projection of an uncertainty set $\mathcal{P}$ to $\mathcal{S}'$ is defined as,
\begin{equation}\label{projection}\operatorname{Proj}_{\mathcal{S}^{\prime}}\mathcal{P}\triangleq \left\{\mathbf{p}_{\mathcal{S}^{\prime}}~\left|~ \exists~ \mathbf{p}_{\mathcal{S}\setminus\mathcal{S}^{\prime}}\colon~(\mathbf{p}_{\mathcal{S}\setminus\mathcal{S}^{\prime}},\mathbf{p}_{\mathcal{S}^{\prime}})\in \mathcal{P}\right\}\right..
    \end{equation}
    Given $\mathbf{p}_{\mathcal{S}^{\prime}} \in \operatorname{Proj}_{\mathcal{S}^{\prime}}\mathcal{P}$ in \eqref{projection}, the conditional projection of  $\mathcal{P}$ to $\mathcal{S}\setminus\mathcal{S}^{\prime}$ with respect to $\mathbf{p}_{\mathcal{S}^{\prime}}$ is defined as 
    \begin{equation}
    \label{conditionalproj}
    \mathcal{P}_{\mathcal{S}\setminus\mathcal{S}^{\prime}}(\mathbf{p}_{\mathcal{S}^{\prime}})\triangleq\{\mathbf{p}_{\mathcal{S}\setminus\mathcal{S}^{\prime}}\mid(\mathbf{p}_{\mathcal{S}\setminus\mathcal{S}^{\prime}},\mathbf{p}_{\mathcal{S}^{\prime}})\in\mathcal{P}\}.
    \end{equation}
    In essence, the conditional projection $\mathcal{P}_{\mathcal{S}\setminus\mathcal{S}^{\prime}}(\mathbf{p}_{\mathcal{S}^{\prime}})$ in \eqref{conditionalproj} contains all possible parameters (i.e., transition probabilities) of state space $\mathcal{S}\setminus\mathcal{S}^{\prime}$ given $\mathbf{p}_{\mathcal{S}^{\prime}}$, and \eqref{projection} ensures $\mathbf{p}_{\mathcal{S}^{\prime}}$ is well-defined.
\end{definition}

From a geometric perspective, the conditional projection can be viewed as the intersection of $\mathcal{P}$ with an affine subspace where $\mathbf{p}_{\mathcal{S}^{\prime}}$  equals the underlying true parameters $\mathbf{p}^0_{\mathcal{S}^{\prime}}$. In other words, the decision-maker can trust the probabilistic information with respect to the state (sub-)space $\mathcal{S}^{\prime}$.

If the uncertainty set $\mathcal{P}$ is rectangular, for any $\mathcal{S}^*\subseteq \mathcal{S}$, we readily know $\mathcal{P}_{\mathcal{S}^*}$ is unrelated to $\mathbf{p}_{\mathcal{S}^{\prime}}$ when $\mathcal{S}^*$ and $\mathcal{S}^{\prime}$ are disjoint, namely, we have $\mathcal{P}=\mathcal{P}_{\mathcal{S}\setminus(\mathcal{S}^{\prime}\cup\mathcal{S}^*)}\times\mathcal{P}_{\mathcal{S}^*}\times \mathcal{P}_{\mathcal{S}^{\prime}}$ by the definition of rectangularity. Conversely, without rectangular assumptions, we cannot obtain such independence. The idea of the  $k$-rectangular ambiguity sets is limiting the number of possible conditional projections of $\mathcal{S}^*$ given $\mathbf{p}_{\mathcal{S}^{\prime}}$.

\begin{definition}[$k$-rectangular ambiguity sets]
    For any $\mathcal{S}^{*}\subseteq \mathcal{S}$, let $\mathcal{S}^{\prime}$ range over all subsets of $\mathcal{S}\setminus\mathcal{S}^{*}$ and $\mathbf{p}_{\mathcal{S}^{\prime}} \in \operatorname{Proj}_{\mathcal{S}^{\prime}}\mathcal{P}$,  a class of conditional projection sets $\mathfrak{P}_{\mathcal{S}^{*}}$ with respect to $\mathcal{S}^*$ is defined as
    \begin{equation}
        \label{k-rectangular}
        \mathfrak{P}_{\mathcal{S}^{*}}\triangleq\left\{\mathcal{P}_{\mathcal{S}^{*}}(\mathbf{p}_{\mathcal{S}^{\prime}})\colon \forall \mathcal{S}^{\prime}\subseteq\mathcal{S}\setminus\mathcal{S}^{*},\forall \mathbf{p}_{\mathcal{S}^{\prime}} \in\operatorname{Proj}_{\mathcal{S}^{\prime}}\mathcal{P}\right\}.
    \end{equation}
    The ambiguity set $\mathcal{P}$ is called $k$-rectangular if $|\mathfrak{P}_{\mathcal{S}^{*}}|\leq k$ for all possible subsets $\mathcal{S}^{*}\subseteq \mathcal{S}$.
\end{definition}

As we discussed above, we can view the given information $\mathbf{p}_{\mathcal{S}^{\prime}}$ as trustful, and the deviations may occur on the space $\mathcal{S}^*$. Limiting the cardinality of the conditional projection sets reduces the possible scenarios. For example, a LDST set constrained by cardinality number $k$ is $(k+1)$-rectangular. Note that (\ref{k-rectangular}) is generalized from (\ref{lightning}), both the probabilistic guarantee and tractability proof adhere to \cite{mannor2012lightning} in a similar way.

Different from restricting the cardinality of deviations, \cite{wiesemann2013robust} resort to \textit{linear decision rules}, an approach adopted in the RO literature \citepalias{ben2009robust}, using constant or (piece-wise) linear functions to approximate the value-to-go function. For a non-rectangular ambiguity set $\mathcal{P}$, the authors first construct $\mathcal{\bar{P}}:=\bigotimes_{s\in \mathcal{S}}\mathcal{P}_s$ as the smallest $s$-rectangular ambiguity set that contains $\mathcal{P}$. Then, they propose an algorithm to optimize the linear approximations of the value-to-go function over the rectangularized ambiguity set and yield an optimal but overly conservative solution in polynomial time.

Motivated by Example \ref{nonrectangularexample}, \cite{gohDataUncertaintyMarkov2018} and \cite{goyalRobustMarkovDecision2023} argue that the transition probability is always affected by environment features. Thanks to big data technology, the authors can develop a new framework to explicitly characterize the relationship between features/factors and transition probability, and they call this modeling framework a \textit{factor uncertainty model}. A factor matrix ambiguity set is defined as
\begin{equation}
    \left.\mathbf{p}=\left\{\left(\textstyle\sum\limits_{i=1}^ru_{sa}^iw_{i,s^{\prime}}\right)_{sas^{\prime}}\right|\boldsymbol{W}=(\boldsymbol{w}_1,\ldots,\boldsymbol{w}_r)\in\mathcal{W}\subseteq\mathbb{R}^{S\times r}\right\},
\end{equation}
where coefficients $u_{1},...,u_{S}$ are fixed and known in $\mathbb{R}_{+}^{r\times \mathcal{A}}$ and underlying factors $\mathcal{W} \in \mathbb{R}_{+}^{\mathcal{S}\times r}$ is a convex, compact subset such that 
\begin{equation}
\textstyle\sum\limits_{i=1}^ru_{sa}^i=1,\forall(s,a)\in\mathcal{S}\times\mathcal{A},~\textstyle\sum\limits_{s'=1}^Sw_{i, s'}=1,\forall i\in[r],
\end{equation}
We would like to highlight that, akin to the approach outlined in \cite{wiesemann2013robust}, this model represents another instance of the utilization of linear decision rules. However, in this context, linear decision rules are applied to approximate the transition probability.

As state-action pairs $(s, a)$ are inefficient to characterize the true dynamics in non-rectangular settings, \citetalias{tirinzoni2018policy} propose a robust learning algorithm with \textit{policy conditioned marginal uncertainty sets}, which incorporates the critical transitions into the construction of uncertainty sets. Specifically, inspired by inverse reinforcement learning (IRL) literature, the authors construct an empirical sample statistics $\hat{\kappa}$ to capture the features of transition tuples $(s, a,s^{\prime})$, instead of directly establishing an empirical distribution as the nominal transition probabilities as same as previous works:
\begin{equation}
	\hat{\kappa} = \frac{1}{N}\sum\limits_{i=1}^N\sum\limits_{t=1}^{T}\phi\left(s_t^{(i)},a_t^{(i)},s_{t+1}^{(i)}\right)
\end{equation}
where $\phi(\cdot):\mathcal{S}\times\mathcal{A}\times\mathcal{S}\rightarrow \mathbb{R}$ is a predetermined feature function. Note that the trajectory samples are generated by a known policy $\tilde{\pi}$, called \textit{base policy}, under the unknown underlying dynamics $\mathbf{p}^0$. The feature function $\phi$ essentially evaluates the ``potential value" of each transition $(s, a,s^{\prime})$ induced by $\tilde{\pi}$ and $\mathbf{p}^0$. In practice, the feature function $\phi$ is often chosen as an indicator function. In this case, a transition $(s, a,s^{\prime})$ evaluated as 1 represents a crucial event, while those evaluated as 0 can be ignored. This setup can be conceptually linked to LDST or $k$-rectangular set. The empirical statistics $\hat{\kappa}$ approximates the feature expectation $\kappa_{\phi}(\tilde{\pi},\mathbf{p}^0):=\mathbb{E}_{\mathbf{p}^0}^{\Tilde{\pi}}\left[\sum_{t=1}^{T-1}{\phi}(S_t,A_t,S_{t+1})\right]$. Given $\hat{\kappa}$, they construct the uncertainty sets as
\begin{equation}
\label{policy_condition_1}
	\mathcal{P}\triangleq\left\{\mathbf{p}\in (\Delta^{|\mathcal{S}|})^{|\mathcal{S}|\times |\mathcal{A}|}~\middle|~ \kappa_{\phi}(\tilde{\pi},\mathbf{p})=\hat{\kappa} \right\},
\end{equation}
or a more general form as
\begin{equation}
\label{policy_condition_2}
	\mathcal{P}\triangleq\left\{\mathbf{p}\in (\Delta^{|\mathcal{S}|})^{|\mathcal{S}|\times |\mathcal{A}|}~\middle|~ \|\kappa_{\phi}(\tilde{\pi},\mathbf{p})-\hat{\kappa}\|\le \theta \right\}.
\end{equation}
The uncertainty sets in \eqref{policy_condition_1} or \eqref{policy_condition_2} require a feasible transition probability matrix $\mathbf{p}$ performs as same as the unknown $\mathbf{p}^0$ in terms of feature expectation with feature function $\phi$ and base policy $\tilde{\pi}$. The change from $(s, a)$ pair to $(s, a, s^{\prime})$ triple essentially involves the comprehensive incorporation of entire trajectories. Finally, the authors demonstrate that a robust control problem with a mixed objective with \textit{policy conditioned marginal uncertainty sets} can be solved in polynomial time. Note that the new issue is not equivalent to the original objective of RMDPs, however, the novel formulation still offers a robust solution in a broad sense.

Recently, \citetalias{wang2023policy} proposed a double-loop algorithm for generic RMDPs, where the outer loop updates the policies and the inner loop updates the worst-case transition
probabilities, alternately. In the outer loop, for a policy $\pi_t$ at iteration $t$, the algorithm finds a transition $\mathbf{p}_t$ of the worst-case one such that
\begin{equation*}
	\mathcal{J}(\pi_t,\mathbf{p}_t)\ge \max_{\mathbf{p}\in \mathcal{P}} \mathcal{J}(\pi_t,\mathbf{p}) -\epsilon_t
\end{equation*}
holds, where $\mathcal{J}(\pi,\mathbf{p}):=\mathbb{E}^{\pi}_{\mathbf{p}}\left[\textstyle \sum\limits _{t=1}^{\infty}\gamma^{t-1}R(s_{t},a_{t},s_{t+1})\right]$ denotes the cumulative discounted reward under the policy $\pi$ and transition matrix $\mathbf{p}$, and $\epsilon_t$ is the predetermined tolerance parameter that satisfies $\epsilon_{t+1}\le \gamma\epsilon_t$. Once $\mathbf{p}_t$ is determined, it is considered as the worst-case transition probabilities, and a projected gradient step is taken with respect to $\pi$ to minimize $\mathcal{J}(\pi,\mathbf{p}_t)$ like in non-robust MDPs, thereby obtaining optimized policy $\pi_{t+1}$. A significant contribution is that the outer loop procedure does not require any rectangular assumption, while the procedure for finding $\mathbf{p}_t$ is not straightforward. Similarly, the inner loop aims to find the worst-case $\mathbf{p}$ given the fixed outer policy $\pi_{t+1}$. However, the inner maximization remains challenging due to its non-convexity, coinciding with previous results that the policy evaluation procedure is computationally difficult \citep{wiesemann2013robust}. To address these challenges, the authors introduce rectangular assumptions specifically for the inner loop in this paper.

\section{Other Related Frameworks and New Fashions in RMDPs}\label{sec7}

What we have introduced so far still falls within the traditional min-max framework, emphasizing the pursuit of an optimal solution to hedge against the (conditional) worst-case. However, the worst-case performance of these approaches is sometimes overly conservative, offering little practical insight into the actual performance of reliable policies. In the RMDPs community, new fashionable modeling approaches have emerged, most inspired by recent advancements in RO and DRO. We will introduce these emerging frameworks and related literature as an extension for future research directions.

\subsection{Satisfactory Regime}
RMDPs can mitigate the impact of perturbations, however, the optimized worst-case performances are often too pessimistic and not easy to offer insights and interpretation for decision-makers \citepalias{long2023robust}. To address this issue, the satisfactory RMDPs have emerged. In the satisfactory regime, there are two mainstream formulations: \textit{Safe Policy Improvement} (SPI) and \textit{Robust Satisficing} (RS).

SPI is established in Batch Reinforcement Learning, where only limited trajectories and a behavioral/baseline policy $\pi_{B}$ to collect the trajectories are accessible. The goal of SPI is obtaining a policy that is guaranteed to perform at least as well as $\pi_{B}$. The general SPI problem can be described as follows.
\begin{definition}[Safe Policy Improvement Problem]
    Given the uncertainty set $\mathcal{P}$ and a baseline policy $\pi_B$, find a maximal $\theta>0$ and a new policy $\pi$ such that $\mathcal{J}(\pi,\mathbf{p})\geq\mathcal{J}(\pi_{B},\mathbf{p})+\theta,~ \forall \mathbf{p}\in\mathcal{P}$, where $\mathcal{J}(\pi,\mathbf{p})$ denotes the cumulative discounted reward under the policy $\pi$ and transition matrix $\mathbf{p}$. The above statement is mathematically equivalent to the following optimization problem:
    \begin{equation}\label{SPI}
\pi_S\in\arg\max_{\pi}\min_{\mathbf{p}\in\mathcal{P}}\left[\mathcal{J}(\pi,\mathbf{p})-\mathcal{J}(\pi_B,\mathbf{p})\right].
    \end{equation}
\end{definition}
Compared to traditional max-min formulation, the problem related to \eqref{SPI} aims to find a policy that maximizes the minimum improvement over the baseline policy $\pi_{B}$, considering all possible realizations in the uncertainty set. It is termed as ``safe" because it guarantees an improvement (as $\theta > 0$) even in the worst-case scenario within the defined uncertainty set.

The seminal paper in SPI by \citetalias{ghavamzadeh2016safe} considers $L_1$-norm uncertainty sets where the discrepancy is bounded by an error function $e(s, a)$ for all $(s, a)\in\mathcal{S}\times\mathcal{A}$:
\begin{equation}
	\mathcal{P}_{sa}\triangleq\left.\left\{\mathbf{p}_{sa}\in\Delta^{|\mathcal{S}|} \right| \|\hat{\mathbf{p}}_{sa}-\mathbf{p}_{sa}\|_1\leq e(s,a) \right\}.
\end{equation}
 The authors prove that the optimal solution to the problem (\ref{SPI}) may be purely randomized by constructing a counterexample, and demonstrate that solving \eqref{SPI} is NP-hard, even though the uncertainty set is $(s, a)$-rectangular. Fortunately, they find that, if the error function induced by the baseline policy is zero, i.e., $e(s,\pi_B(s))=0$ for all $s\in\mathcal{S}$, the original problem \eqref{SPI} can be simplified and solved in polynomial time. Motivated by this finding, the authors propose a simple and practical approximate algorithm in which the error function is updated as
\begin{equation*}
	\tilde{e}(s,a)=\begin{cases}
	e(s,a), &\text{when }\pi_\mathrm{B}(s)\neq a\\
	0, &\mathrm{otherwise}
	\end{cases}
	\quad \forall ~(s,a)\in\mathcal{S}\times\mathcal{A}.
\end{equation*}

Subsequently, \citetalias{laroche2019safe} and \citetalias{simao2020safe} both generalize SPI as SPI with baseline bootstrapping, effectively mitigating the unreliability caused by the large state space and insufficient samples. The key procedure to achieve this is distinguishing the types of state-action pairs. Specifically, when a state-action pair $(s, a)$ is rarely seen in the dataset (indicating high uncertainty), the trained policy $\pi$ will default to the baseline policy $\pi_B$ such that $\pi(a|s)=\pi_B(a|s)$. One can see it as a \textit{knows-what-it-knows} algorithm, seeking assistance from the baseline policy when uncertain, instead of searching for the analytic optimum over the whole space. Particularly, while \cite{laroche2019safe} assume the baseline policy $\pi_B$ is known, \cite{simao2020safe} further relax this assumption and just assume there is an MLE estimate for $\pi_B$. Leveraging the statistical properties of bootstrapping, \cite{laroche2019safe} propose an efficient algorithm that provides PAC-style guarantees of policy improvement with high probability, and \cite{simao2020safe} offer SPI guarantees over the true baseline $\pi_B$ even without direct access to it.

Robust Satisficing (RS) is a novel framework developed in RO and DRO \citep{long2023robust} that transforms the traditional optimization problem of seeking the optimum into attaining an acceptable target specified ahead. Due to the outstanding interpretability and solid theoretical foundations of RS, \citetalias{ruan2023robust} extend it into the MDPs version recently. 

Recall that solving an MDP can be formulated as an LP. To deliver the core idea of RS, we compare the robust version of dual LP and RS version of dual LP, omitting the corresponding primal form for simplicity. The robust version of dual LP can be computed as
\begin{equation}
\label{robustLP}
	 \begin{array}{ll}
	 \max & \boldsymbol{r}^\top\boldsymbol{u}\\
	 \mathrm{s.t.} & \sum\limits_{a\in\mathcal{A}} u_{sa}-\gamma\sum\limits_{s'\in\mathcal{S}}\sum\limits_{a\in\mathcal{A}}p_{s'as}u_{s'a}-d_s\leq0  ~\forall~\mathbf{p}\in\mathcal{P},s\in\mathcal{S}\\\end{array}
\end{equation}
where $\boldsymbol{r}=(R(s,a))_{s\in\mathcal{S},a\in\mathcal{A}}$ are the immediate rewards and $u_{sa}$ is the dual variable with respect to state-action pair $(s,a)$. Here, uncertainty set $\mathcal{P}=\{\mathbf{p}\in (\Delta^{|\mathcal{S}|})^{|\mathcal{S}|\times |\mathcal{A}|} \mid \ell(\mathbf{p}, \hat{\mathbf{p}})\le \theta\}$. As opposed to traditional RMDPs with hard constraints to limit the distinction among distributions, a target-oriented RSMDP imposes soft constraints for all other transition distributions:
\begin{equation}
\label{rsLP}
	 \begin{array}{ll}
	 	\min & \sum_{s\in\mathcal{S}} k_s\\
	 	\mathrm{s.t.} & \sum\limits_{a\in\mathcal{A}} u_{sa}-\gamma\sum\limits_{s'\in\mathcal{S}}\sum\limits_{a\in\mathcal{A}}p_{s'as}u_{s'a}-d_s\leq k_s\cdot \ell(\mathbf{p}, \hat{\mathbf{p}})  ~\forall~s\in\mathcal{S}\\
	 	& \sum\limits_{s\in\mathcal{S}}p_{sa}=1 ~\forall~s\in\mathcal{S}\\
	 	& \boldsymbol{r}^\top\boldsymbol{u}\ge \tau
	 	\end{array}
\end{equation}
where $\tau$ is the pre-specified target. To interpret the meanings of the decision variables $k_s$ in \eqref{rsLP}, we introduce the concept of 
constraint violation, where the constraints in \eqref{robustLP} do not hold. Hence, the decision variables $\{k_s\}$ reflect the magnitude of constraint violation incurred by the discrepancy between nominal and true transition probabilities. For instance, if the values of $\{k_s\}$ are small enough, then $k_s\cdot \ell(\mathbf{p}, \hat{\mathbf{p}})$ will also tend to be small such that the only mild violation will occur. However, as the existence of $\tau$, the feasible region of $\{k_s\}$ is limited, which in turn to restricts the discrepancy $\ell(\mathbf{p}, \hat{\mathbf{p}})$.

Compared to \eqref{robustLP} which tries to attain the maximal ``rewards", the RSMDP attempts to minimize the magnitude of the constraint violation while ensuring that the “reward” remains at least $\tau$. This approach underscores the concept of satisficing, prioritizing acceptable outcomes over purely optimal rewards. Leveraging the theorems from \cite{long2023robust}, the authors formulate the original problems as a conic program when the distance function $\ell$ is chosen as $\ell_1$-norm. Furthermore, inspired by the FOM framework proposed by \cite{grand2021scalable}, they transform \eqref{rsLP} into a min-max form and apply a PDA to solve the min-max problem efficiently.

\subsection{Online and Bayesian Regime}

As the focus of MDPs or RMDPs is \textit{planning}, another significant drawback of traditional RMDPs lies in their inability to ``learn" from uncertainty, which easily causes overly conservative policies due to the extensive uncertainty sets. To tackle this challenge and infuse the adaptive learning characteristic of RL, there is a growing trend to incorporate online learning or Bayesian concepts into the design of uncertainty/ambiguity sets, Bellman updates, environmental configurations, and other facets of RMDPs. Such efforts are expanding the scope of RMDPs and injecting fresh dynamism into the field.

Out of consideration for limiting the size of uncertainty set, \citetalias{lim2013reinforcement} assume there is a unique but unknown subset $\mathcal{Z}$ of state-action pairs behaving adversarially while all others are stochastic, as in non-robust MDPs. In essence, they propose a two-step framework: divide state-action pairs into two categories first and then handle them respectively. Motivated by \cite{bubeck2012best} that investigate a similar scenario in a multi-armed bandit setting, the authors employ a statistical hypothesis test, termed \textit{consistency check}, after executing arbitrary state-action pair $(s, a)$ every time. All state-action pairs that pass the consistency check are stochastic, which means the incurred rewards and transitions are safely accepted, while the pairs that do not pass the check will be included in $\mathcal{Z}$. Combining the consistency check with UCRL2, a well-known online algorithm that treats every state-action pair as stochastic \citep{auer2009near}, the authors develop two algorithms for finite-horizon and infinite-horizon cases, respectively. Following the optimistic property of UCRL2, the results of the proposed algorithms are less conservative than the standard RMDPs.

Sharing similar concerns but adopting distinct methodologies, \cite{petrik2019beyond} observe that reliance on concentration inequalities results in conservative uncertainty sets. To circumvent this issue, 
the authors propose leveraging Bayesian inference to construct tighter uncertainty sets. They focus on a safe return measured by Value-at-Risk (VaR) and aim to improve the lower bound when the value function is unknown. Under Bayesian assumptions, the transition probabilities $\mathbf{p}_{sa}$ is a random variable with a prior distribution. Specifically, they introduce the \textit{robustification with sensible value functions} (RSVF) concept, in which the uncertainty set $\mathcal{P}_{sa}$ is derived from a set of possible value functions $\mathcal{V}$. For any given value function $v\in \mathcal{V}$, they approximate the optimal uncertainty set $\mathcal{K}_{sa}(v)$ as
\begin{equation}
    \mathcal{K}_{sa}(v)=\left\{\mathbf{p}_{sa}\in \Delta^{|\mathcal{S}|}~\left|~\mathbf{p}_{sa}^{\top}v\leq \operatorname{VaR}_{\hat{\mathbf{p}}}^{\theta}[\hat{\mathbf{p}}_{sa}^{0\top}v]\right\}\right.,
\end{equation}
where $\hat{\mathbf{p}}$ serves as the nominal probability. This set essentially defines a set of plausible transition probabilities $\mathbf{p}_{sa}$ whose expected value meets the VaR criterion for being ``sufficiently safe'' regarding the current value function $v$. 

Note that $\mathcal{K}_{sa}(v)$ is sufficient to provide a reasonable uncertainty set if the value function $v$ is known. However, in the context of this paper, the value function is unknown. As new value functions are progressively added to $\mathcal{V}$, RSVF updates the uncertainty set $\mathcal{P}_{sa}$ to append plausible $\mathbf{p}_{sa}$ that is ``close'' enough to the worst-case distribution within the current $\mathcal{K}_{sa}(v), v\in\mathcal{V}$. When $\mathcal{P}_{sa}$ is updated, we can derive the (currently) optimal value function $\hat{v}^*$ and policy $\hat{\pi}^*$. If the intersection $\mathcal{K}_{sa}(\hat{v}^*)\cap\mathcal{P}_{sa}\neq \emptyset$ for all $(s, a)$ pairs, which means that $\hat{v}^*$ is ``sufficiently robust'' for some adversarial distributions, the algorithm deems the current policy $\hat{\pi}^*$ to be safe. They provide a finite iteration probabilistic guarantee and empirically compare the safe estimates from RSVF and other non-Bayesian methods. Empirical results demonstrate that RSVF indeed provides much tighter ambiguity sets and safe returns.

Unlike the aforementioned works, \citetalias{derman2020bayesian} offer a posterior perspective to design Bellman update for Bayesian RMDPs. Collecting an observation history $\mathcal{T}:=\{(s_1,a_1),...,(s_t,a_t)\}$ induced by policy $\pi$ and a predetermined confidence level $\theta_{sa}>0$ for each pair, the posterior uncertainty sets are constructed over time:
\begin{equation}
    \widehat{\mathcal{P}}^{t}_{sa}\triangleq \left\{\mathbf{p}_{sa}\in \Delta^{|\mathcal{S}|} \big|  \|\mathbf{p}_{sa}-\bar{\mathbf{p}}_{sa}\|_{1}\le \theta_{sa} \right\},
\end{equation}
where $\bar{\mathbf{p}}_{sa}:=\mathbb{E}[\mathbf{p}_{sa}|\mathcal{T}]$ acts as the nominal transition probability, and the overall ambiguity set satisfies $(s, a)$-rectangular assumption. Once the posterior uncertainty set is established, a posterior over robust Q-values at stage $t$ can then be computed as:
\begin{equation}
    \label{URBE}
    \widehat{Q}^{t}(s,a)=R(s,a)+\gamma\inf_{\mathbf{p}\in\widehat{\mathcal{P}}^{t}_{sa}}\textstyle\sum\limits_{s^{\prime} \in \mathcal{S}, a^{\prime} \in \mathcal{A}}\pi^t(a^{\prime}\mid s^{\prime})p_{sas^{\prime}}\widehat{Q}^{t+1}(s^{\prime},a^{\prime}).
\end{equation}
Equation (\ref{URBE}) is termed as \textit{Uncertainty Robust Bellman Equation} (URBE). This uncertainty set is constructed for each episode based on the observed data from all previous ones.

Slightly modifying the procedures in \citetalias{o2018uncertainty}, the authors compute posterior variance $\omega_{sa}^t$ of $\widehat{Q}^{t}(s, a)$ in (\ref{URBE}), and approximate the posterior over robust Q-values as normal distribution $\mathcal{N}(\bar{Q}, \mathrm{diag}(\boldsymbol{\omega}))$, where $\mathrm{diag}(\boldsymbol{\omega})$ is the solution of URBE and $\bar{Q}$ the conditional expectation of $\hat{Q}$. This approach offers a trade-off between robustness and conservatism for robust policies. Besides, the authors propose a DQN-URBE algorithm, and they show that it is significantly faster to change dynamics online compared to existing robust techniques with fixed uncertainty sets.

Contrary to most studies that concentrate on developing innovative robust algorithms to learn from the original environment, \citetalias{wang2023robust} introduce an RL framework designed to approximate an adversarial environment and generate adversarial data for training purposes. This approach is inspired by \citetalias{kumar2024policy} which reveals that the adversarial distribution fundamentally alters the next-state transition probability from the nominal one. Compared to applying intricate iterations of robust methods, directly employing non-robust methods in adversarial environments will be more computationally efficient, thus addressing scalability challenges in RMDPs with high-dimensional domains. Different from the $L_p$-norm ambiguity sets used by \cite{kumar2024policy}, the authors construct KL-divergence ambiguity sets under $(s, a)$-rectangular assumption, which circumvents the issue of different supports between nominal and adversarial distribution. They also derive a closed-form expression for the adversarial distribution, depicted as a re-weighted nominal distribution. Consequently, the framework can be straightforwardly implemented by simulating the next state under the nominal dynamic and choosing it with corresponding adversarial weight. After realizing the next state, $(s, a,s')$ tuple is added to the data buffer, and the policy is trained with data from the buffer via any non-robust RL method.

\subsection{Risk, Regularization, and Robustness}

Our pursuit of robustness aims to enhance performance and mitigate uncertainty when the underlying truth is ambiguous. This concept resonates across various fields, e.g., finance and machine learning (ML),  where both risk-aware and regularized formulations have seen considerable success. Notably, a body of literature attempts to bridge these concepts with RO and DRO \citep{shapiro2017distributionally,namkoong2017variance,shapiro2021tutorial,shafieezadeh2019regularization,blanchet2019quantifying}, and also RMDPs. These studies offer diverse interpretations of robustness, revealing connections that shed light on the models of MDPs. This understanding helps to partially alleviate the conservatism inherent in robust policies. 

In terms of risk-aware MDPs, \cite{osogami2012robustness} explores the equivalence between RMDPs and risk-aware MDPs with an expected exponential utility objective. By leveraging the properties of this objective, the author reformulates risk-aware MDPs into a general RMDP that aims to minimize the expected cumulative cost adjusted by the KL-divergence penalty. This penalty quantifies the discrepancy between the nominal and the worst-case distribution. Building on this foundation,  \cite{bauerleDistributionallyRobustMarkov2022} broaden these equivalence findings to encompass all risk-aware MDPs that employ spectral risk measures, which are known as a class of coherent risk measures like the \textit{Expected Shortfall}. However, the risk measures examined are confined to the Markov risk measure \citep{ruszczynski2010risk}, known to be potentially non-gradient-dominant \citepalias{huang2021convergence}. 

More recently, \citetalias{zhang2023regularized} establish the equivalence between risk-aware MDPs and a class of regularized RMDPs, including the standard RMDPs. The authors introduce an innovative formulation for risk-aware MDPs that incorporates general convex risk measures, requiring only that these measures adhere to principles of monotonicity, translation invariance, and convexity. This formulation facilitates a broader equivalence and supports the use of policy gradient methods for algorithm design, ensuring global convergence.

Similarly, the equivalence between regularization and robustness in RMDPs has been investigated recently. \cite{derman2020distributional} consider a reward-maximizing RMDP within the Wasserstein distance ambiguity set. Employing strong duality and necessary conditions, they construct a regularized value function that serves as a lower bound to the distributionally robust value. Moreover, they extend the finite-sample guarantee results of \cite{yang2020wasserstein} to regularized MDPs, though they stop analyzing the tightness and asymptotic consistency of their approach as the sample size increases.

As a supplement and improvement of \cite{derman2020distributional}, \citetalias{derman2021twice} formally demonstrate the equivalence between regularized MDPs and RMDPs with arbitrary ball-constrained ambiguity sets, extending beyond the Wasserstein metric. They propose a comprehensive RMDP model that accounts for uncertainties in both rewards and transitions, introducing a novel extension to regularized MDPs that encompasses both policy and value regularization. Utilizing contraction and monotonicity assumptions, they successfully apply Banach’s fixed point theorem to this extended model and show that equivalence is maintained. As a byproduct, they propose efficient Bellman updates for a modified policy iteration tailored to the structure of the regularized problem. While their method solves RMDPs similarly to classical non-robust MDPs, the obtained policies can be overly conservative due to the lack of restrictions on worst-case transition probabilities to the probability simplex.

Because of theoretical equivalence and computational efficiency, the integration of risk measures and regularization into MDPs has recently gained significant momentum. Value-at-Risk (VaR) and Conditional Value-at-Risk (CVaR), as the most widely recognized risk metrics, have spurred a considerable amount of research in quantile-based MDPs. In the model-based domain, \citetalias{li2022quantile} investigate quantile MDP(QMDP) which is defined as
\begin{equation}\label{QMDP}
    \max_{\pi\in\Pi}Q_\tau^\pi\left[\sum_{t=0}^{T-1}r_t(s_t,a_t,\xi_t)\right],
\end{equation}
where $Q_\tau(X):=\inf\{x\mid P(X\leq x)\geq\tau\}$ and $\tau\in(0,1)$ is the quantile value. The authors reformulate (\ref{QMDP}) as a max-min optimization and have proven that the quantile-based value-to-go function is equivalent to solving a specific optimization problem for each action. This breakthrough facilitates the application of dynamic programming techniques. Additionally, they craft efficient algorithms and conduct a detailed complexity analysis. In terms of model-free approaches, 
\cite{yu2022risk} adapt CVaR as a special one-step conditional risk measure to preserve time consistency and subsequently reformulate the risk-averse MDP as a risk-neutral counterpart with augmented action space and the adjusted immediate rewards. Furthermore, the authors show that the Bellman operator in this formulation is a contraction mapping, and consequently, the authors extend Q-learning results into this case and develop a risk-averse deep Q-learning algorithm.

The concept of weaving regularization into dynamic programming also has been thoroughly investigated within the reinforcement learning field. \citep{kaufman2013robust, neu2017unified,kostrikov2021offline,ho2021partial, kumar2022efficient}.  While these regularized variants of value iteration and policy iteration have been developed, \cite{grand2022convex} is the first to propose the convex optimization formulation of RMDPs. This innovative approach integrates an entropic regularization of the robust Bellman operator with a change of variables involving exponential and logarithmic functions. By conducting this convex reformulation, arbitrary algorithms designed for solving convex optimization problems can be incorporated with the reformulated RMDPs. This novel perspective has paved the way for the development of innovative algorithms tailored to RMDPs.

\section{Conclusion}
This survey provides an overview of the theory, methodologies, and advances of RMDPs with uncertain transition probability. Before the details of the literature on RMDPs, we elaborate on the foundations of RMDPs and DRMDPs and introduce an important concept \textit{rectangularity}. Rectangularity plays a key role in most RMDPs research, allowing decomposability to ensure that RMDPs can be solved in polynomial time. We summarize three class definitions of rectangularity and propose a new and straightforward proof that solving non-rectangular RMDPs is NP-hard. 

Mainly, we categorize the ambiguity modeling approaches into three groups, parametric, moment-based, and discrepancy-based. The parametric approach is the simplest yet most limited way among the three methods to characterize the uncertainty. Since it usually imposes individual support constraints or limits the distribution family, it exhibits statistically poor performance and heavily relies on LPs to compute the optimal policy. Due to the issues of interpretability, moment-based ambiguity modeling is more suitable for situations where the randomness is exogenous or concrete problem backgrounds. As the increasingly popular approach, discrepancy-based modeling enjoys desirable statistical guarantees and performance in a data-driven manner. We summarize three prevalent discrepancy measures--norm, $\phi$-divergence, and Wasserstein distance. All of them provide powerful finite sample guarantees to contain the underlying probability measure within constructed ambiguity sets with (high) probability, in turn, guiding how many samples we need to make a safe decision. It should not be ignored that the research of fast algorithms based on discrepancy-based ambiguity sets has gained more and more attention.

Finally, we review recent efforts that depart from rectangular assumptions and minimax framework. Relaxing traditional rectangular assumptions allows capturing of dependencies between transitions, which better characterizes the reasonable impact of uncertainty and reduces conservatism. However, these coupled uncertainty models often need to compromise flexibility and generality. Besides the minimax objective, there are also other approaches to hedge uncertainty and achieve robust performance. These new paradigms bring new techniques and perspectives, complementing the research of RMDPs. 

Sequential decision-making under uncertainty is a difficult yet crucial research problem in practice. Although RMDPs have a substantial theoretical basis, their applications in sequential decision-making are relatively limited compared to other tools. One potential explanation for this limitation is rectangularity. How to reasonably relax this assumption and develop new frameworks is an important direction, particularly within a data-driven context and when integrating artificial intelligence techniques. Meanwhile, the development of fast algorithms and the provably robust RL algorithms will be an interesting and practical direction.

% References here (outcomment the appropriate case)

% CASE 1: BiBTeX used to constantly update the references
%   (while the paper is being written).
%\bibliographystyle{informs2014} % outcomment this and next line in Case 1
%\bibliography{<your bib file(s)>} % if more than one, comma separated

% CASE 2: BiBTeX used to generate mypaper.bbl (to be further fine tuned)
%\input{mypaper.bbl} % outcomment this line in Case 2

%If you don't use BiBTex, you can manually itemize references as shown below.
\newpage
\bibliographystyle{agsm}
\bibliography{bibfile}

%%%%%%%%%%%%%%%%%
\end{document}